\documentclass[12pt]{article}
\usepackage{amsmath, amscd, amssymb, latexsym,hhline}
\usepackage[all]{xy}
\def\:{{\colon}}

\def\FPdim{{\mbox{\rm FPdim\,}}}
\def\coev{{\mbox{\rm coev}}}
\def\ev{{\mbox{\rm ev}}}
\def\tr{{\mbox{tr}}}
\def\ker{{\mbox{ker\,}}}
\def\Tr{{\mbox{Tr}}}

\def\mod{{\bf mod }}

\def\BC{{\mathbb C}}

\def\Irr{{\mbox{\rm Irr}}}

\def\a{{\alpha}}
\def\b{{\beta}}
\def\d{{\delta}}

\def\g{{\gamma}}

\def\H'{{D^{\w'}(G)}}
\def\H{{D^\w(G)}}

\def\e{{\varepsilon}}

\def\w{{\omega}}

\def\map{{\longrightarrow}}
\def\qed{{\ \ \ $\square$}}
\def\pf{{\it Proof\ \ }}

\def\<{{\langle}}
\def\>{{\rangle}}
\def\pf{{\it Proof. }}
\def\C{{\,\subseteq\,}}

\newcommand{\ldual}[1]{{^*\!#1}}

\newtheorem{thm}{Theorem}[section]
\newtheorem{cor}[thm]{Corollary}
\newtheorem{defn}[thm]{Definition}
\newtheorem{prop}[thm]{Proposition}
\newtheorem{remark}[thm]{Remark}

\newtheorem{lem}[thm]{Lemma}

\newtheorem{example}[thm]{Example}

\numberwithin{equation}{section}
\newcommand{\vspc}{\vspace{.15 in}}

\begin{document}
\title{Central Invariants and Frobenius-Schur Indicators for Semisimple Quasi-Hopf Algebras}
\author{Geoffrey Mason and Siu-Hung Ng}
\date{}
\maketitle
\markright{  \hspace{5cm}            \sc G. Mason and
S.-H. Ng}{}
\begin{abstract}
In this paper, we obtain a canonical central element $\nu_H$ for
each semi-simple quasi-Hopf algebra $H$ over any field $k$ and prove that $\nu_H$ is invariant under
gauge transformations. We show that if $k$ is algebraically closed of characteristic zero then
for any irreducible representation of $H$ which affords
the character $\chi$, $\chi(\nu_H)$ takes only the values 0, 1 or -1,
moreover if $H$ is a Hopf algebra or a twisted quantum double of a
finite group then $\chi(\nu_H)$ is the corresponding Frobenius-Schur
Indicator.  We also prove
an analog of a Theorem of Larson-Radford for split semi-simple quasi-Hopf algebra over any field $k$.
Using this result, we establish the relationship between the antipode $S$,
the values of $\chi(\nu_H)$, and certain associated bilinear forms when the underlying field $k$ is
algebraically closed of characteristic zero.

\end{abstract}
\section{Introduction}
   In the paper (\cite{LM00}), Linchenko and Montgomery introduced
and studied Frobenius-Schur
indicators for irreducible representations of a semi-simple Hopf
algebra $H$ over an algebraically closed
field of characteristic $p \neq 2$. If
  $\Lambda$ is the \emph{unique} normalized left integral of $H$, i.e.
$\e(\Lambda) = 1$, set
\begin{equation}\label{eq: centralhopfelt}
\nu = \nu_H = \sum_{(\Lambda)} \Lambda_1 \Lambda_2.
\end{equation}
Here we have used Sweedler notation $  \Delta(\Lambda) =
\sum_{(\Lambda)} \Lambda_1 \otimes \Lambda_2$, so that
if $m$ is multiplication in $H$ then $\nu = m \circ \Delta
(\Lambda)$. Then $\nu$ is a central element of $H$ and
the Frobenius-Schur indicator $\nu_{\chi}$ of an irreducible
$H$-module $M$ with character $\chi$ is defined via
\begin{equation}\label{eq: hopfs}
\nu_{\chi} = \chi(\nu).
\end{equation}
In case $H$ is a group algebra $k[G],  \nu = |G|^{-1}\sum_{g \in
G} g^2$ and $\nu_{\chi} = |G|^{-1}\sum_{g \in G} \chi(g^2)$
reduces to the original definition of Frobenius and Schur
(cf.\cite{CR88} or \cite{serre77}, for example). Generalizing the
famous result of Frobenius and Schur for group algebras, Linchenko
and Montgomery show that for general semi-simple $H$, $\nu_{\chi}$
can take only the values $0, 1$, or $-1$. Moreover $\nu_{\chi}
\neq 0$ if, and only if, $M \cong M^*$, and in this case $M$
admits a non-degenerate $H$-invariant bilinear form $\<\cdot
,\cdot \>$ satisfying
\begin{equation}
  \<u, v\> = \nu_{\chi} \<v, u\>
\end{equation}
for $u, v \in M$. Recall that $\<\cdot ,\cdot \>$ is $H$-invariant if
\begin{equation}\label{eq: invbilform}
  \sum_{(h)}\<h_1u, h_2v\> = \epsilon(h)\<u, v\>
\end{equation}
for $h \in H$ and $u, v \in M$.

\vspc

\noindent In a recent paper (\cite{KMM02}) the authors showed how
one may effectively compute Frobenius-Schur indicators for a
certain class of Hopf algebras. Their work applies, in particular,
to the case of the \emph{quantum double} $D(G)$ of a finite group
$G$, and it was shown (loc. cit.) how the indicators for
irreducible modules over $D(G)$ may be given in terms of purely
group-theoretic invariants associated to $G$ and its subgroups.
The algebra $D(G)$ is of interest in orbifold conformal field
theory (\cite{Mas}), indeed in this context there is a more
general object, the \emph{twisted quantum double} $D^\omega(G)$,
that arises naturally (\cite{DPR90}). (Here, $\omega \in Z^3(G,
\BC^\times)$ is a normalized $3$-cocycle about which we shall have
more to say below.) The present work originated with a natural
problem: understand Frobenius-Schur indicators for twisted quantum
doubles.

\vspc
\noindent $D^{\omega}(G)$ is a semi-simple \emph{quasi-Hopf} algebra
(over $\BC$, say), but is generally
not a Hopf algebra. One of the
  difficulties this imposes is that the antipode $S$ is not
necessarily involutorial (something that is always true
for semi-simple Hopf algebras by a Theorem of Larson and Radford (\cite{LaRa87})), whereas
having $S^2 = id$ is fundamental for the
Linchenko-Montgomery approach and therefore for the calculations in
(\cite{KMM02}. If it happens that $S^2 = id$
  then Theorem 4.4 of (loc. cit.) can be used to obtain indicators given by
\begin{equation} \label{eq: kmmfsind}
\nu_{\chi} = |G|^{-1} \sum_{ x^{-1}gx=g^{-1}} \gamma_x(g, g^{-1})
\theta_g(x, x) \chi(e(g) \otimes x^2).
\end{equation}
(Undefined notation is explained below; $\gamma_x$ and $\theta_g$ are
certain $2$-cochains determined by $\omega$.)

\vspc
\noindent If $G$ is \emph{abelian} then $D^{\omega}(G)$ is a Hopf
algebra (\cite{MN01}), though perhaps with a non-trivial
$\beta$ element, and for any
$G$ it turns out that one can always \emph{gauge} $\omega$, i.e.
replace it by a cohomologous $3$-cocycle $\omega'$, in such a way
that the antipode for $D^{\omega'}(G)$ is an involution. So (\ref{eq:
kmmfsind}) provides a preliminary solution to our problem,
but it is unsatisfactory for the following reason: if we gauge $\omega$,
the new $3$-cocycle $\omega'$ will give new values for the
Frobenius-Schur indicators which in general are
not the same as the original values. While this may not be an issue
if one is interested in a fixed $D^{\omega}(G)$,
there are both mathematical and physical reasons for insisting that
the FS indicators for $D^{\omega}(G)$ be \emph{robust},
that is they depend only on the cohomology class of $\omega$. From
this standpoint, (\ref{eq: kmmfsind}) is generally not what
we are looking for. We need a more \emph{functorial} approach.

\vspc
\noindent One knows that if $\omega$ and $\omega'$ are cohomologous
then $D^{\omega}(G)$ and
$D^{\omega'}(G)$ are gauge-equivalent and that therefore the
corresponding module categories are tensor equivalent (cf. \cite{Drin90},
\cite{DPR90}, \cite{Kassel}).
Indeed, it follows from a result of Etinghof and Gelaki (\cite{EG02}) that
the converse is also true, so that gauge-equivalence of the twisted
doubles is the \emph{same} as tensor equivalence of the module categories. So
we are looking for invariants of such module categories with respect
to tensor equivalence.
Because Hopf algebras and twisted doubles are not closed with respect
to gauge equivalence,
  this means that we have to work with the module categories
of
\emph{arbitrary} semi-simple quasi-Hopf algebras.

\vspc
\noindent Peter Bantay has introduced a notion of  indicator into
rational conformal field theory from a rather different
point-of-view (\cite{Bantay97}, \cite{Bantay00}). His point of departure is the Verlinde formula
  and the $S$ and $T$ matrices associated to a RCFT. To this modular
data together with an irreducible character
he associates a certain
numerical expression and shows that it is equal once again to either
$0, 1$ or $-1$. It is
possible to evaluate Bantay's indicator in case the  matrices $S$ and
$T$ are associated to a twisted double
$D^{\omega}(G)$ (\cite{BantayConv})
and one obtains the expression
\begin{equation}\label{eq: bantayfsind}
  |G|^{-1} \sum_{ x^{-1}gx=g^{-1}}\omega(g^{-1}, g, g^{-1})
\gamma_x(g, g^{-1}) \theta_g(x, x) \chi(e(g) \otimes x^2).
\end{equation}
Compared to (\ref {eq: kmmfsind}), (\ref{eq: bantayfsind}) contains
an extra term $\omega(g^{-1}, g, g^{-1})$. Furthermore,
it is easy to see that (\ref{eq: bantayfsind}) is robust in the previous sense.

\vspc
\noindent Suppose that $H$ is any semi-simple quasi-Hopf algebra, and
let $M$ be an irreducible $H$-module
with character $\chi$. In the present paper we will
   construct a \emph{canonical} central element $\nu_H$ of $H$ with
the following properties:
\begin{quote}
  \begin{enumerate}
\item[(a)] $ \nu_H$ is invariant under \emph{any} gauge
transformation of  $H$\,. \item[(b)]If $H$ is a Hopf algebra then
$\nu_H$ coincides with (\ref{eq: centralhopfelt})\,.
\item[(c)]$
\mbox{If} \ H = D^{\omega}(G) \ \mbox{then} \ \chi(\nu_H)  \
\mbox{coincides with Bantay's indicator (\ref{eq: bantayfsind})}.$
\item[(d)] Asumme that $k$ is algebraically closed and char $k=0$.
\begin{enumerate}
\item[(i)] $\chi(\nu_H) = 0, 1, \ \mbox{or} -1  $.
\item[(ii)] $  \chi(\nu_H) \neq 0$ if, and only if,  ${^*\!M} \cong M$\,.
In this case, $M$ admits a certain  non-degenerate   bilinear form
$\<\cdot , \cdot\> $ such that
\begin{equation}\label{qhopfbilform}
\<x, y\> = \<y, g^{-1} x\>
\end{equation}
for all $x$, $y \in M$. Here, $g$ is a distinguished element of $H$, which we call the
\emph{trace element}, which is independent of $M$.
\item[(iii)] $\displaystyle \Tr(S) =
\sum_{\chi \in \Irr(H)}\chi(\nu_H) \chi(\b^{-1})$\,.
\end{enumerate}
\end{enumerate}
\end{quote}

Part (d) is the analog for general semi-simple quasi-Hopf algebras
of the corresponding result in \cite{LM00} for Hopf algebras. The bilinear
form $\<\cdot , \cdot\> $
has a certain adjointness property with respect to the antipode $S$ of
$H$, and there are relations to an analog of a Theorem of
Larson-Radford ($S$ is involutorial for semi-simple Hopf algebras).
Namely, we show that for a semi-simple quasi-Hopf algebra the
antipode is involutorial up to conjugation.
The trace element $g$ plays an important role in our discussion of (d),
in particular its properties lead to the fact that the category
$H$-$\mod_{fin}$ of finite-dimensional H-modules is a pivotal category in
the sense of Joyal and Street. For twisted doubles, $g$ coincides with
$\beta$, while for Hopf algebras the
Larson-Radford Theorem implies that $g = 1$.\\

  The proof that  $\chi(\nu_H)$ takes only the values 0, 1 or -1 is
somewhat elaborate. Indeed, in an earlier version of the present
paper (\cite{MN03}) this had been left open.
Subsequently, Pavel Etingof alerted us to the existence of his recent
preprint with Nikshych and Ostrik \cite{ENO} on fusion categories, and
suggested that some of the results obtained there could be used
to help settle the issue of the values of our indicator. More
precisely, Etingof pointed out that our trace element $g$ defines an
isomorphism of tensor functors $Id \longrightarrow {^{**}?}$. This together with
$S(g) =g^{-1}$ are the main ingredients in the proof.\\

  The paper is organized as follows: we cover some basic facts about
quasi-Hopf algebras in Section 2, including several strategically
important elements in $H \otimes H$ introduced by Hausser and Nill
\cite{HNQA}. In Sections 3 and 4 we define the central element $\nu_H$
and establish that the family of
Frobenius-Schur indicators $\chi(\nu_H)$ is a gauge invariant for
semi-simple quasi-Hopf algebras. In section 5 we show that our
indicators coincide with those of Bantay in the case of a twisted
double. In Section 6 we introduce the trace element $g$ and establish
the analog of the Larson-Radford Theorem, while Section 7 is devoted
to further properties of $g$ as discussed above. Section 8 covers the
relation of indicators to bilinear forms and completes the proof of
(d)(i), and in  Section 9 we return to the case of twisted doubles to
complete the analysis in that case. For simplicity, we will only work on
algebraically fields of characteristic zero in Section 7, 8 and 9.\\

The authors are indebted to Pavel Etingof for his interest and
extended correspondence, and also thank Peter Bantay and Susan
Montgomery for helpful discussions.

\section{Quasi-Hopf Algebras}
In this section we recall the definition of quasi-Hopf algebras
and their properties described in \cite{Drin90} and \cite{Kassel}.
Moreover, we recall some interesting results recently obtained in
\cite{HNQA}, \cite{HN991},\cite{HN992} and \cite{PV00}. In the
sequel, we will use the notation introduced in this section.
Throughout this paper, we will always assume that $k$ is a field
and any algebras and vector spaces are over $k$. In section 7, 8 and 9, we will further assume
$k$ to be an algebraically closed field of characteristic zero.\\

A {\em quasi-bialgebra} over $k$ is a 4-tuple $(H, \Delta, \e,
\Phi)$, in which $H$ is an algebra over $k$, $\Delta: H \map H
\otimes H $ and $\e : H \map k$ are algebra maps, and $\Phi$ is an
invertible element in $H \otimes H \otimes H$ satisfying the
following conditions:
\begin{gather}
   \label{1.1}
  (\e \otimes id)\Delta(h) =h= (id \otimes \e)\Delta(h);  \\
   \label{1.2}
  \Phi (\Delta \otimes id)\Delta(h)\Phi^{-1} = (id \otimes
  \Delta)\Delta(h) \mbox{ for all } h \in H; \\
   \label{1.3}
  (id \otimes  id  \otimes \Delta)(\Phi)
  (\Delta \otimes id \otimes  id )(\Phi) =
  (1 \otimes \Phi)( id \otimes\Delta \otimes  id )(\Phi)(\Phi
  \otimes 1);\\
   \label{1.4}
  (id \otimes \e \otimes id)(\Phi) = 1 \otimes 1\,.
\end{gather}
The maps $\Delta$, $\e$ and $\Phi$ are respectively called the diagonal
map, counit, and associator of the quasi-bialgebra. If there is no
ambiguity, we will simply write $H$ for the quasi-bialgebra $(H,
\Delta, \e, \Phi)$. Using (\ref{1.3}), one can also easily see
that
\begin{gather}
   \label{1.5}
  (\e \otimes id \otimes id)(\Phi) = 1 \otimes 1 =
  (id \otimes id \otimes \e)(\Phi) \,.
\end{gather}
Moreover, the module category $H$-\mod of the quasi-bialgebra $H$
is a tensor category (cf. \cite{Drin90} and \cite{Kassel} for the
details). \\

Following \cite{Kassel}, a {\em gauge transformation} on a
quasi-bialgebra $H=(H, \Delta, \e, \Phi)$ is an invertible element
$F$ of $H \otimes H$ such that
$$
(\e \otimes id)(F)=1 =(id \otimes \e)(F)\,.
$$
Using a gauge transformation on $H$, one can define an algebra map
$\Delta_F: H \map H \otimes H$ by
\begin{equation}\label{2.6}
  \Delta_F(h) = F \Delta(h) F^{-1}
\end{equation}
for any $h \in H$, and an invertible element $\Phi_F$ of $H \otimes
H \otimes H$ by
\begin{equation}\label{2.7}
  \Phi_F=(1 \otimes F) (id \otimes \Delta)(F) \Phi (\Delta \otimes
  id)(F^{-1})(F^{-1} \otimes 1)\,.
\end{equation}
Then $H_F=(H, \Delta_F, \e, \Phi_F)$ is  also a quasi-bialgebra. \\

Two quasi-bialgebras $A$ and $B$ are said to be {\em gauge
equivalent} if there exists a gauge transformation $F$ on $B$ such
that $A$ and $B_F$ are isomorphic as quasi-bialgebras. If $A$ and
$B$ are gauge equivalent quasi-bialgebras,  $A$-$\mod$,
$B$-$\mod$ are equivalent tensor categories (cf. \cite{Kassel}). Conversely, if
$A$, $B$ are finite-dimensional semi-simple quasi-bialgebra such that
$A$-\mod and $B$-\mod are equivalent tensor categories, then $A$ and $B$ are
gauge equivalent quasi-bialgebras (cf. \cite{EG02}).\\

A quasi-bialgebra  $( H, \Delta, \e,\Phi)$ is called a {\em
quasi-Hopf algebra} if there exist an anti-algebra automorphism
$S$ of $H$ and elements $\a, \b \in H$ such that for all element
$h \in H$, we have
\begin{equation}\label{2.8}
  \sum_{(h)} S(h_1)\a h_2 =\e(h)\a, \quad \sum_{(h)} h_1\b S(h_2)
  =\e(h)\b\,
\end{equation}
and
\begin{equation}\label{2.9}
  \sum_i X_i \b S(Y_i)\a Z_i =1, \quad
  \sum_i S(\overline{X}_i) \a \overline{Y}_i\b S(\overline{Z}_i) =1
\end{equation}
where $\Phi=\sum_i X_i \otimes Y_i\otimes Z_i$, $\Phi^{-1}= \sum_i
\overline{X}_i \otimes \overline{Y}_i\otimes \overline{Z}_i$ and
$\sum_{(h)} h_1\otimes  h_2 = \Delta(h)$\,. We shall write $( H,
\Delta, \e,\Phi, \a, \b ,  S)$ for the complete data of the
quasi-Hopf algebra and $S$ is called the \emph{antipode} of $H$. When the
context is clear, we will simply write $H$ for the quasi-Hopf
algebra $(H, \Delta, \e,\Phi, \a, \b , S)$. One can easily see
that a Hopf algebra is a quasi-Hopf algebra with $\Phi=1 \otimes 1
\otimes 1$ and $\a=\b=1$.\\

Unlike a Hopf algebra, the antipode for a quasi-Hopf algebra is generally not unique.
\begin{prop}{\rm \cite[Proposition 1.1]{Drin90}}\label{prop2.1}
Let $H=(H, \Delta, \e,\Phi, \a, \b , S)$ be a quasi-Hopf algebra.
If $u$ is a unit of $H$ then $H_u= (H, \Delta, \e,\Phi, u \a, \b u^{-1} , S_u)$
is also a quasi-Hopf algebra, where $S_u(h)=uS(h)u^{-1}$ for all $h \in H$.
Conversely, for any $\a', \b' \in H$ and for any algebra anti-automorphism $S'$
of $H$ such that $H'=(H, \Delta, \e,\Phi, \a', \b' , S')$ is a quasi-Hopf algebra,
then there exist a unique invertible element $u$ of $H$ such that
$$
H_u = H'\,.
$$
\hfill\qed
\end{prop}

If $F$ is a gauge transformation on the quasi-Hopf algebra $H=(H,
\Delta, \e,\Phi, \a, \b ,  S)$, we can define $\a_F$ and $\b_F$ by
$$
\a_F = \sum_i S(d_i) \a e_i
\quad\mbox{and}\quad \b_F = \sum_i f_i \b S(g_i)
$$
where $F=\sum_i f_i \otimes g_i$  and $F^{-1} =\sum_i
d_i \otimes e_i$. Then, $H_F=(H, \Delta_F,
\e,\Phi_F, \a_F, \b_F ,  S)$ is also a quasi-Hopf algebra.\\

The antipode of a Hopf algebra is known to be a anti-coalgebra map.
For a quasi-Hopf algebra $H$, this is true up to conjugation.
Following \cite{Drin90}, we define $\g, \d \in H \otimes H$ by the formulae
\begin{eqnarray}
\g & = & \sum_i S(U_i)\a V_i \otimes S(T_i) \a W_i\,,\\
\d & = & \sum_j K_j \b S(N_j) \otimes L_j \b S(M_j)\,,
\end{eqnarray}
where
\begin{eqnarray*}
\sum_i T_i \otimes U_i \otimes V_i \otimes W_i &=& (1 \otimes \Phi^{-1})(id
\otimes id \otimes \Delta)(\Phi) \,,\\
\sum_j K_j \otimes L_j \otimes M_j \otimes N_j &=& (\Delta \otimes  id
\otimes id )(\Phi)( \Phi^{-1}\otimes 1 ) \,.
\end{eqnarray*}
Then,
\begin{equation}\label{FH}
F_H=\sum_i (S \otimes S)(\Delta^{op}(\overline{X}_i)) \cdot \g \cdot
\Delta(\overline{Y}_i \b S(\overline{Z}_i))
\end{equation}
is an invertible element of $H \otimes H$ where $\Phi^{-1}= \sum_i
\overline{X}_i \otimes \overline{Y}_i\otimes \overline{Z}_i$.
Moreover,
$$
F_H \Delta(S(h))F_H^{-1} = (S \otimes S)(\Delta^{op}(h)\,.
$$
for all $h \in H$.\\

The category of  finite-dimensional left $H$-module of a
quasi-Hopf algebra $H$ with antipode $S$, denote by
$H$-$\mod_{fin}$ is a rigid tensor category. Let $M$ be a
finite-dimensional left $H$-module and $M'$ its $k$-linear dual. Then
the $H$-action on $M'$, given by
$$
(h\cdot f)(m) = f(S(h)m)
$$
for any $f \in M'$ and $m \in M$, defines a left $H$-module
structure on $M'$. We shall denote  by
$\ldual{M}$ the left dual of the $M$ in $H$-$\mod_{fin}$.
Similarly, the right dual of $M$, denote by $M^*$, is the
$H$-module with the underlying $k$-linear space $M'$ with the
$H$-action  given by
$$
(h\cdot f)(m) = f(S^{-1}(h)m)
$$
for any $f \in M'$ and $m \in M$ (cf. \cite{Drin90}).\\

In \cite{HNQA}, \cite{HN991} and \cite{HN992}, Frank Hausser and
Florian Nill introduced some interesting elements in $H \otimes H$
for any arbitrary quasi-Hopf algebra $H=(H, \Delta, \e,\Phi, \a,
\b , S)$ in the course of studying the corresponding theories of
quantum double, integral and the fundamental theorem for
quasi-Hopf algebras. These elements of $H \otimes H$ are given by
\begin{gather}
\label{R}
q_R  = \sum X_i \otimes S^{-1} (\a Z_i)Y_i \,, \quad p_R= \sum
\overline{X}_i \otimes \overline{Y}_i \b
S(\overline{Z}_i)\,, \\
\label{L}
q_L = \sum S(\overline{X}_i) \a \overline{Y}_i \otimes
\overline{Z}_i\,, \quad  p_L = \sum Y_iS^{-1}(X_i\b) \otimes Z_i
\,
\end{gather}
where $\Phi= \sum_i X_i \otimes Y_i \otimes Z_i$ and $\Phi^{-1}= \sum_i
\overline{X}_i \otimes \overline{Y}_i\otimes \overline{Z}_i$. One can show
easily (cf. \cite{HNQA}) that  they obey the relations (for all $a
\in H$)
\begin{eqnarray}
 \label{eq:qR}
 (a\otimes 1)\, q_R &= &  \sum (1 \otimes S^{-1}(a_2))\,q_R\, \Delta(a_1)  , \\
\label{eq:qL}
(1 \otimes a)\, q_L &= & \sum (S(a_1)\otimes 1)\, q_L \, \Delta(a_2)  , \\
 \label{eq:pR}
 p_R\, (a \otimes 1) &= &  \sum \Delta(a_1)\, p_R \, (1 \otimes  S(a_2)) , \\
  \label{eq:pL}
p_L\, (1 \otimes a) &= &   \sum\Delta(a_2)\,p_L\, (S^{-1}(a_1)\otimes
1)\,.
\end{eqnarray}
where $\Delta(a)= \sum a_1 \otimes a_2$.
Suppressing the summation symbol
and indices, we write $q_R = q^1_R \otimes q^2_R$, etc.
These elements also satisfy the identities (cf. \cite{HNQA}):
\begin{eqnarray}
\label{2.16}
\Delta(q_R^1) \, p_R\, (1\otimes S(q_R^2)) &= & 1\otimes 1, \\
\label{2.17}
(1\otimes S^{-1}(p_R^2))\, q_R \, \Delta(p_R^1)& = & 1\otimes 1, \\
\label{2.18}
\Delta (q_L^2) \, p_L \, (S^{-1}(q_L^1) \otimes 1)&= &1 \otimes 1, \\
\label{2.19}
(S(p_L^1) \otimes 1)\, q_L \, \Delta(p_L^2) & = & 1\otimes 1.
\end{eqnarray}
We will use these equations in the sequel.
\section{Central Gauge Invariants for Semi-simple Quasi-Hopf Algebras}
Suppose that $H=(H, \Delta, \e, \Phi, \a, \b, S)$ is a
finite-dimensional quasi-Hopf algebra. A left integral  of $H$ is
an element $l$ of $H$ such that $h l = \e(h) l$ for all $h \in H$.
A right integral of $H$ can be defined similarly. It follows from
\cite{HNQA}  that the subspace of left (right) integrals  of $H$
is of dimension 1. Moreover, if $H$ is semi-simple, the subspace
of left integral is identical to the space of right integrals of
$H$ and $\e(\Lambda) \ne 0$ for any non-zero left
integral $\Lambda$ of $H$ (see also  \cite{PV00}). We will call the
two-sided integral $\Lambda$ of $H$ {\em normalized} if $\e(\Lambda)=1$.\\

Let $\Lambda$ be a left integral of $H$. Then for any $a \in H$,
\begin{equation}\label{eq:left_integral}
\e(a)\Delta(\Lambda) =
\Delta(a) \Delta(\Lambda)\,.
\end{equation}
Similarly, if $\Lambda'$ is a right integral of $H$, then we have
\begin{equation}\label{eq:right_integral}
\e(a)\Delta(\Lambda') =
 \Delta(\Lambda') \Delta(a)\,.
\end{equation}

We then have the following lemma.
\begin{lem}\label{pq_1}
Let $H=(H, \Delta, \e, \Phi, \a, \b, S)$ be a finite-dimensional quasi-Hopf
algebra.
\begin{enumerate}
\item[\rm (i)] If $\Lambda$ is a left integral of $H$, then for any $a \in H$,
\begin{eqnarray}
\label{3.3}
(1 \otimes a)q_R \Delta(\Lambda) &=& (S(a) \otimes 1)q_R \Delta(\Lambda)\,,\\
\label{3.4}
(1 \otimes a)q_L \Delta(\Lambda) & =& (S(a) \otimes 1)q_L \Delta(\Lambda)\,, \\
\label{3.5}
\mbox{and}\quad (\b \otimes 1)q_L \Delta(\Lambda)&=& (\b \otimes 1)q_R \Delta(\Lambda)=\Delta(\Lambda)\,.
\end{eqnarray}
\item[\rm (ii)] If $\Lambda'$ is a right integral of $H$, then for any $a \in H$,
\begin{eqnarray}
\label{3.6}
\Delta(\Lambda')p_R (a \otimes 1) & = &   \Delta(\Lambda')p_R (1 \otimes S(a))
,\\
\label{3.7}
\Delta(\Lambda')p_L (a \otimes 1)  &= & \Delta(\Lambda')p_L (1 \otimes S(a))\,,\\
\label{3.8}
\mbox{and}\quad\Delta(\Lambda') p_L (1 \otimes \a)&=&\Delta(\Lambda') p_R
(1 \otimes \a)=\Delta(\Lambda')\,.
\end{eqnarray}
\end{enumerate}
\end{lem}
\pf (i) By the equations (\ref{eq:left_integral}) and (\ref{eq:qR}), for any
$a \in H$,
\begin{eqnarray*}
(a\otimes 1)q_R \Delta(\Lambda) & = &  (1 \otimes S^{-1}(a_1)))q_R
\Delta(a_2)\Delta(\Lambda)\\
&=&(1 \otimes S^{-1}(a_1\e(a_2)))q_R
\Delta(\Lambda)\\
&=& ( 1 \otimes S^{-1}(a))q_R\Delta(\Lambda) .
\end{eqnarray*}
Hence, by substituting $a$ with $S(a)$, we prove equation (\ref{3.3}).
Now we have
\begin{eqnarray*}
\Delta(\Lambda) &=& (1\otimes S^{-1}(p_R^2))\, q_R \, \Delta(p_R^1)\Delta(\Lambda)\quad \mbox{ (by (\ref{2.17}))}\\
&=& (1\otimes S^{-1}(p_R^2\e(p_R^1))q_R \Delta(\Lambda) \quad \mbox{ (by (\ref{eq:left_integral}))}\\
& = & (1 \otimes S^{-1}(\b))q_R \Delta(\Lambda) \quad \mbox{ (by (\ref{1.5}))}\\
&=&(\b \otimes 1)q_R \Delta(\Lambda) \quad \mbox{ (by (\ref{3.3}))}\,.
\end{eqnarray*}
The remaining formulae in (i) and (ii) can be proved similarly using equations
(\ref{1.5}), (\ref{eq:qR})-(\ref{2.19}),  (\ref{eq:left_integral}) and
(\ref{eq:right_integral}). \qed\\

\begin{lem}\label{l3.2}\label{invariant_1}
Let $H=(H, \Delta, \e, \Phi, \a, \b, S)$ be a finite-dimensional quasi-Hopf
algebra and $F$ a gauge transformation on $H$.
Suppose that $q_R^F, q_L^F, p_R^F, p_L^F$ are  the corresponding
$p$'s and $q$'s for $H_F$ defined in {\rm (\ref{R})} and {\rm (\ref{L})}.
\begin{enumerate}
\item[\rm (i)] If $\Lambda$ is a left integral of $H$, then
$$q_R^F \Delta_F(\Lambda)=q_R\Delta(\Lambda)F^{-1} \,\quad
\mbox{and} \quad
q_L^F \Delta_F(\Lambda)=q_L\Delta(\Lambda)F^{-1}\,.
$$
\item[\rm (ii)] If $\Lambda'$ is a right integral of $H$, then
$$
 \Delta_F(\Lambda')p_R^F =  F \Delta(\Lambda')p_R
\,\quad
\mbox{and} \quad
\Delta_F(\Lambda') p_L^F =  F \Delta(\Lambda')p_L\,.
$$
\end{enumerate}
\end{lem}
\pf (i) Let $\Phi^{-1}=\sum_j \overline{X}_j \otimes \overline{Y}_j \otimes
\overline{Z}_j$, $F=\sum_i f_i \otimes g_i$ and $F^{-1}=\sum_l d_l \otimes e_l$.
Then, we obtain
\begin{eqnarray*}
\Phi_F^{-1} & = &(F \otimes 1)(\Delta  \otimes id)(F) \Phi^{-1}
(id \otimes \Delta)(F^{-1})(1 \otimes F^{-1}) \\
& =  &(F \otimes 1)\left(\sum_{i,j,l} f_{i,1} \overline{X}_j d_l \otimes
                        f_{i,2} \overline{Y}_j e_{l,1} \otimes
                        g_i \overline{Z}_j e_{l,2}\right) (1 \otimes F^{-1})
\end{eqnarray*}
where $\Delta(f_i)=\sum f_{i,1} \otimes f_{i,2}$ and
$\Delta(e_l)=\sum e_{l,1} \otimes e_{l,2}$. Thus, we have
\begin{eqnarray*}
q_L^F \Delta_F(\Lambda) & = & \left(\sum S(f_{i'}f_{i,1} \overline{X}_jd_l)\a_F g_{i'}f_{i,2} \overline{Y}_je_{l,1}
\otimes g_i \overline{Z}_j e_{l,2}\right)F^{-1} F\Delta(\Lambda) F^{-1} \\ \\
      & = & \left(\sum S(f_{i,1} \overline{X}_jd_l)\a f_{i,2} \overline{Y}_je_{l,1}
\otimes g_i \overline{Z}_j e_{l,2}\right)\Delta(\Lambda) F^{-1}
\quad(\mbox{since } \sum S(f_{i'}) \a_F g_{i'} = \a)\\ \\
      & = & \left(\sum S(f_{i,1} \overline{X}_jd_l \e(e_l))\a f_{i,2} \overline{Y}_j
\otimes g_i \overline{Z}_j \right)\Delta(\Lambda) F^{-1} \quad (\mbox{by }(\ref{eq:left_integral}))\\ \\
& = & \left(\sum S(f_{i,1} \overline{X}_j)\a f_{i,2} \overline{Y}_j
\otimes g_i \overline{Z}_j \right)\Delta(\Lambda) F^{-1}
\quad(\mbox{since } \sum d_l\e(e_l) = 1_H)\\
&=& \left(\sum S(\overline{X}_j)\a \e(f_i) \overline{Y}_j
\otimes g_i \overline{Z}_j \right)\Delta(\Lambda) F^{-1}
\quad(\mbox{since } \sum S(f_{i,1})\a f_{i,2} = \e(f_i)\a) \\ \\
&=& \left(\sum S(\overline{X}_j)\a \overline{Y}_j
\otimes \overline{Z}_j \right)\Delta(\Lambda) F^{-1}
\quad(\mbox{since } \sum \e(f_i)g_i=1_H )\\ \\
& =&  q_L \Delta(\Lambda) F^{-1})
\end{eqnarray*}
The other three equations can be proved similarly. \qed \\

\begin{thm} \label{th3.3}Let $H=(H, \Delta, \e, \Phi, \a, \b, S)$
be a finite-dimensional quasi-Hopf algebra. Suppose that $\Lambda$ is a two-sided integral
of $H$. Then, the elements
\begin{equation*}\label{eq:invariant_1}
q_R \Delta(\Lambda) p_R\,,  \quad  \quad q_R \Delta(\Lambda) p_L\,,\quad
q_L \Delta(\Lambda) p_R\,, \mbox{ and }  \quad q_L \Delta(\Lambda) p_L\
\end{equation*}
in $H \otimes H$ are invariant under gauge transformations. Moreover,
$$
m(q_R \Delta(\Lambda) p_R)=m(q_R \Delta(\Lambda) p_L)=m(q_L \Delta(\Lambda) p_R)=m(q_L \Delta(\Lambda) p_L)
$$
where $m$ denote the multiplication of $H$. In addition, $m(q_R \Delta(\Lambda) p_R)$ is a central
element of $H$.
\end{thm}
\pf It follows from Lemma \ref{invariant_1} that for any gauge transformation
 $F$ on $H$,
$$
q_*^F\Delta_F(\Lambda)=q_*\Delta(\Lambda)F^{-1},\quad\mbox{and} \quad
\Delta_F(\Lambda)p_*^F = F\Delta(\Lambda)p_*\,
$$
where $q_*^F =q_L^F$ or $q_R^F$ and $p_*^F =p_L^F$ or $p_R^F$.
Thus we have
\begin{eqnarray*}
q_*^F\Delta_F(\Lambda)p_*^F &= &q_*\Delta(\Lambda)F^{-1}p_*^F \\
& = & q_*F^{-1}\Delta_F(\Lambda)p_*^F \\
&= & q_* F^{-1}F\Delta(\Lambda)p_*\\
&= & q_*\Delta(\Lambda)p_*\,.
\end{eqnarray*}

Let $m$ denote the multiplication of $H$ and let $m_{RR}$, $m_{RL}$, $m_{LR}$ and $m_{LL}$ denote the elements
$$
m(q_R \Delta(\Lambda) p_R), \, m(q_R \Delta(\Lambda) p_L), \, m(q_L \Delta(\Lambda) p_R), \,\,\mbox{ and }
m(q_L \Delta(\Lambda) p_L)
$$
respectively.   Then for any $a \in H$,
\begin{eqnarray*}
S(a)m_{RR}  & = &
m( (S(a) \otimes 1)q_R  \Delta(\Lambda)p_R) \\
& =& m((1\otimes a)q_R \Delta(\Lambda)p_R) \quad \mbox{by Lemma \ref{pq_1}(i)}\\
& =& m( q_R \Delta(\Lambda)p_R(a\otimes 1))\\
& =& m(q_R \Delta(\Lambda)p_R (1\otimes S(a))) \quad \mbox{by  Lemma \ref{pq_1}(ii)}\\
 &=& m_{RR}S(a)\,.
\end{eqnarray*}
As $S$ is an automorphism, the above equation implies that $m_{RR}$ is
in the center of $H$. Using the same kind of arguments, one can
show
that $m_{RL}, m_{LR}$ and $m_{LL}$ are each in the center of $H$.\\

 Let $Q_R$, $Q_L$, $P_R$ and $P_L$ denote the elements
$$
m(q_R \Delta(\Lambda)), \, m(q_L \Delta(\Lambda)),\,
 m(\Delta(\Lambda) p_R)\,, \,\mbox{ and } m(\Delta(\Lambda) p_L)\,
$$
respectively.
Then, we have
\begin{equation}\label{eq:exchange_1}
\begin{aligned}
m_{RR} & = \sum q_R^1 \Lambda_1 p_R^1 q_R^2 \Lambda_2 p_R^2\\
&=\sum S(p_R^1) q_R^1 \Lambda_1  q_R^2 \Lambda_2 q_R^2 \quad \mbox{by Lemma \ref{pq_1}(i)}\\
& = S(p_R^1) Q_R p_R^2 \\
&= \sum_j S(\overline{X}_j) Q_R \b\overline{Y}_jS(\overline{Z}_j)
\end{aligned}
\end{equation}
 and
\begin{equation}\label{eq:exchange_1.2}
\begin{aligned}
m_{RR} & = \sum q_R^1 \Lambda_1 p_R^1 q_R^2 \Lambda_2 p_R^2\\
&=\sum q_R^1 \Lambda_1 p_R^1 \Lambda_2 p_R^2 S(q_R^2) \quad \mbox{by Lemma \ref{pq_1}(ii)}\\
& = q_R^1 P_R S(q_R^2)\\
&= \sum X_i P_R S(Y_i) \a Z_i
\end{aligned}
\end{equation}
where $\Phi^{-1}=\sum_j \overline{X}_j \otimes \overline{Y}_j \otimes
\overline{Z}_j$ and $\Phi=\sum_i X_i \otimes Y_i \otimes
Z_i$.
 Similarly,
\begin{equation}\label{eq:exchange_2}
 m_{LL} = S(\overline{X}_j) \a \overline{Y}_j P_L S(\overline{Z}_j) \\
= \sum X_i \b S(Y_i) Q_L Z_i\,.
\end{equation}
By (\ref{3.5}) and (\ref{3.8}), we have
\begin{equation}\label{eq:PQ}
\begin{aligned}
Q_R = m_{RL}\a, &\quad Q_L =m_{LR}\a \,, \\
P_L = \b m_{RL}, &\quad P_R = \b m_{LR}\,.
\end{aligned}
\end{equation}

Therefore, using equation (\ref{2.9}), we have
$$
m_{RR} = S(\overline{X}_j) m_{RL}\a \overline{Y}_j \b S(\overline{Z}_j) =
m_{RL}S(\overline{X}_j) \a \overline{Y}_j \b S(\overline{Z}_j) = m_{RL}
$$
Similarly, using equations (\ref{eq:exchange_1.2}), (\ref{eq:exchange_2}) and (\ref{eq:PQ})
we can prove
$$
m_{RR}=m_{LR}=m_{LL}\,.
$$
\hfill\qed\\

In \cite{HNQA} and \cite{PV00}, it is shown that a finite-dimensional
quasi-Hopf algebra $H$ is semi-simple if, and only, if there exist a unique
normalized two-sided integral. In this case, we have the following:
\begin{defn}
{\rm Let $H=(H, \Delta, \e, \Phi, \a, \b, S)$ be a
finite-dimensional semi-simple quasi-Hopf algebra and let $\Lambda$ be
the unique normalized two-sided integral of $H$.
We denote by $\nu_H$  the central element
$$
m(q_L \Delta(\Lambda)p_L)
$$
discussed in Theorem \ref{th3.3}.}
\end{defn}

\begin{cor}\label{cor3.5}
Let $H=(H, \Delta, \e, \Phi, \a, \b, S)$ be a finite-dimensional semi-simple quasi-Hopf algebra and
$\Lambda$  the normalized two-sided integral of $H$. Then $\nu_H$
is invariant under gauge transformations, that is
$$
\nu_H = \nu_{H_F}\,
$$
for any gauge transformation $F$ on $H$. Moreover,
$$
\b\a\nu_H = \nu_H\b\a=\sum(\Lambda_1\Lambda_2)
$$
where $\sum \Lambda_1 \otimes \Lambda_2 = \Delta(\Lambda)$. In particular,
if both $\a$ and $\b$ are units of $H$,
then
$$
\nu_H=\sum(\Lambda_1\Lambda_2)(\b\a)^{-1} = (\b\a)^{-1}\sum(\Lambda_1\Lambda_2)\,.
$$
\end{cor}
\pf The first statement follows immediately from Theorem
\ref{3.3}. By equation (\ref{3.5}) and (\ref{3.8}), we have
$\b\nu_H\a= \sum(\Lambda_1\Lambda_2)$. Since $\nu_H$ is central, then the
result
follows. \qed\\

\begin{cor}
Let $H=(H, \Delta,\e, \Phi, \a, \b, S)$,
$H'=(H', \Delta',\e', \Phi', \a', \b', S')$  be  semisimple quasi-Hopf algebras.
If $H$ and $H'$ are
gauge equivalent quasi-bialgebras via the gauge transformation $F$ on $H$ and the
quasi-bialgebra isomorphism $\sigma : H_F \map H'$, then
$$
\sigma(\nu_{H}) = \nu_{H'}\,.
$$
In particular, if $u$ is a unit of $H$, then $\nu_{H_u} = \nu_{H}$.
\end{cor}
\pf Since $H_F$ and $H'$ are isomorphic quasi-bialgebras,
    $(H', \Delta', \e', \Phi', \sigma(\a_F),\sigma(\b_F), \sigma S \sigma^{-1})$ is a
    quasi-Hopf algebra. By Proposition \ref{prop2.1}, there exists a unit $u$ of $H'$
    such that
    \begin{equation}\label{eq:3.13}
    \sigma S \sigma^{-1}(a) = u S'(a) u^{-1}, \quad \sigma(\b_F) = u \a'\, \mbox{ and }
    \sigma(\b_F) = \b'u^{-1}\,.
    \end{equation}
    for all $a \in H'$. Then, we have
 \begin{equation}\label{eq:3.14}
    \sigma S^{-1} \sigma^{-1}(a) = S'^{-1}(u) S'^{-1}(a) S'^{-1}(u^{-1})\,.
  \end{equation}
    Let $\Lambda$ be the normalized two-sided integral of $H$. Since $\sigma$
    is a quasi-bialgebra isomorphism, $\sigma(\Lambda)$ is then a two-sided integral of $H'$ and
    $$
    \e'(\sigma(\Lambda)) =  \e(\Lambda) = 1\,.
    $$
    Therefore, $\Lambda'=\sigma(\Lambda)$ is the unique normalized integral of $H'$. In particular,
    we have
    $$
    (\sigma\otimes \sigma)\Delta_F(\Lambda)= \sum \Lambda'_1 \otimes \Lambda'_2 \quad \mbox{ and }\quad
     (\sigma \otimes \sigma \otimes \sigma )(\Phi_F)= \Phi'
    $$
    where $\sum \Lambda'_1 \otimes \Lambda'_2 = \Delta'(\Lambda')$.
    Let
    $$
     \begin{aligned}
    & \Phi_F = \sum X^F_i \otimes Y^F_i \otimes Z^F_i \,,\quad
    \Phi_F^{-1} =
    \sum \overline{X}^F_j \otimes \overline{Y}^F_j  \otimes \overline{Z}^F_j\,,\\
    & \Phi'= \sum X'_i \otimes Y'_i \otimes Z'_i\,,\quad
    \Phi'^{-1}  =
    \sum \overline{X}'_j \otimes \overline{Y}'_j  \otimes \overline{Z}'_j\,,\\
    &\mbox{and} \quad(\sigma\otimes \sigma)\Delta_F(\Lambda) = \sum \Lambda_1^F \otimes \Lambda_2^F\,.
    \end{aligned}
     $$
    Then,
    $$
    \begin{aligned}
    & \sigma(\nu_H) = \sigma(\nu_{H_F})\quad \quad \mbox{by Corollary (\ref{cor3.5})} \\
    = \,\,& \sigma\left(\sum X^F_i \Lambda^F_1 \overline{X}^F_j S^{-1}(\a_F Z^F_i)
    Y^F_i \Lambda^F_2 \overline{Y}^F_j \b_F S(\overline{Z}^F_j)\right)\\
    = \,\,& \sum X'_i \Lambda'_1 \overline{X}'_j\, (\sigma S^{-1})(\a_F Z^F_i)
    Y'_i \Lambda'_2 \overline{Y}'_j \sigma(\b_F)\, (\sigma S)(\overline{Z}^F_j)\\
    = \,\,& \sum X'_i \Lambda'_1 \overline{X}'_j \,(\sigma S^{-1} \sigma^{-1})(Z'_i) \,
    (\sigma S^{-1}\sigma^{-1})(\sigma(\a_F))
    Y'_i \Lambda'_2 \overline{Y}'_j \sigma(\b_F) \,(\sigma S \sigma^{-1})(\overline{Z}'_j)\\
%   =\,\,& \sum X'_i \Lambda'_1 \overline{X}'_j \,S'^{-1}(u) S'^{-1} (Z'_i) S'^{-1}(u^{-1}) \,
%    S'^{-1}(u)S'^{-1}(u \a'))S'^{-1}(u^{-1})\,
%     Y'_i \Lambda'_2 \overline{Y}'_j \b' u^{-1} \,uS'(\overline{Z}'_j)u^{-1}\\
    =\,\,& \sum X'_i \Lambda'_1 \overline{X}'_j \,S'^{-1}(u) S'^{-1} (Z'_i) S'^{-1}(\a')\,
    Y'_i \Lambda'_2 \overline{Y}'_j \b'S'(\overline{Z}'_j)u^{-1} \quad \mbox{by (\ref{eq:3.13}) and (\ref{eq:3.14})}\\
    = \,& \sum X'_i \Lambda'_1 \overline{X}'_j S'^{-1}(\a' Z'_i)
    Y'_i \Lambda'_2 \overline{Y}'_j \b'S'(\overline{Z}'_j)u u^{-1}  \quad \mbox{by Lemma \ref{pq_1}(ii)}\\
    = \, & \nu_{H'}\,.
    \end{aligned}
    $$

    For any unit $u$ of $H$, $H$ and $H_u$ are obviously gauge equivalent
    as quasi-bialgebra under the gauge transformation $1 \otimes 1$ and the
    quasi-bialgebra isomorphism $id_H$. Hence, the second statement follows.
    \hfill\qed
\section{Frobenius-Schur Indicators}
Let $(H, \Delta,\e, \Phi, \a, \b, S)$ be a semi-simple quasi-Hopf
algebra over the field $k$. Let $M$ be an irreducible $H$-module
with character $\chi$. We call $\chi(\nu_H)$ the Frobenius-Schur
indicator of $\chi$ (or $M$). The family of Frobenius-Schur
indicators $\{ \chi(\nu_H) \}$ is in fact  an invariant of the
tensor category $H$-\mod for any semi-simple quasi-Hopf algebra
$H$.

\begin{thm}\label{thm4.1}
Let $H=(H, \Delta,\e, \Phi, \a, \b, S)$ and $H'=(H', \Delta',\e',
\Phi', \a', \b', S')$ be finite-dimensional semi-simple quasi-Hopf
algebras over an algebraically closed field $k$ of characteristic zero.
If $H$-\mod and $H'$-\mod are equivalent as $k$-linear
tensor categories, then  the families of Frobenius-Schur
indicators for $H$ and $H'$ are identical.
\end{thm}
\pf If $H$-\mod and $H'$-\mod are equivalent as $k$-linear tensor
categories, then, by \cite[Theorem 6.1]{EG02}, $H$ and $H'$ are
gauge equivalent quasi-bialgebras. Suppose that $F$ is a gauge
transformation on $H$ and $\sigma : H_F \map H'$ is a
quasi-bialgebra isomorphism. It follows from Corollary
\ref{cor3.5} that
$$
\sigma(\nu_H)=\nu_{H'}\,.
$$
Let $\Irr(H)$, $\Irr(H')$ be the set of irreducible characters of
$H$ and $H'$ respectively. Then, the map $\chi' \mapsto \chi'
\circ \sigma$ is a bijection from $\Irr(H')$ onto $\Irr(H)$.
Moreover, for any irreducible character $\chi'$ of $H'$,
$$
\chi' \circ \sigma(\nu_H)= \chi'(\nu_{H'})\,.
$$
Thus, $\{\chi'(\nu_{H'}) \}_{\chi' \in \Irr(H')}$ is identical of
the family $\{\chi(\nu_{H}) \}_{\chi \in \Irr(H)}$. \qed

\begin{remark}
{\rm If $H$ is a semi-simple Hopf algebra, then $\Phi=1 \otimes
1\otimes 1$ and $\a=\b=1$. It follows from Corollary \ref{cor3.5}
that
$$
\nu_H= \sum \Lambda_1 \Lambda_2
$$
where $\sum \Lambda_1 \otimes \Lambda_2=\Delta(\Lambda)$ and $\Lambda$ is the normalized
two-sided integral of $H$. Thus, $\chi(\nu_H)$ coincides
with the Frobenius-Schur indicator defined in \cite{LM00}.}
\end{remark}

As an application of Theorem, we give a simple alternative proof
of the fact that $\BC [Q_8]$-\mod and $\BC [D_8]$-\mod are not equivalent as
$\BC$-linear tensor categories where $Q_8$ and $D_8$ are the
quaternion group and the dihedral group of order 8 respectively
(cf. \cite{TaYa98}).

\begin{prop}{\rm \cite{TaYa98}}
The $\BC$-linear categories $\BC[Q_8]$-\mod and $\BC[D_8]$-\mod are not
equivalent as tensor categories.
\end{prop}
\pf Let $G=D_8$ or $Q_8$. Then, $G$ has four degree 1 characters
and one degree 2 irreducible character $\chi_2$. Let $z$ be the
non-trivial central element of $G$. Then $\chi_2(z)=-2$ and
$\chi(z) = 1$ for any character $\chi$ of $G$ of degree 1. Since
$\nu_G=\frac{1}{8}\sum_{g \in G} g^2$, one can easily obtain that
$$
\nu_{Q_8} = \frac{1}{8}(6z+2e), \quad \mbox{and}\quad \nu_{D_8}=
\frac{1}{8}(2z+6e)
$$
where $e$ is the identity of the group. Thus, the family of
Frobenius-Schur indicators for $Q_8$ is $\{1,1,1, 1, -1\}$ but the
the family of Frobenius-Schur indicators for $D_8$ is
$\{1,1,1,1,1\}$. By virtue of Theorem \ref{thm4.1}, $\BC[Q_8]$-\mod
and $\BC[D_8]$-\mod are not equivalent as $\BC$-linear tensor
categories. \qed
\section{Bantay's Formula for Indicators of Twisted Quantum Doubles}
\label{s:Bantay_Indicator}
In this section, we will show that if $H$ is a twisted quantum
double of a finite group $G$ over the field $k$ such that $|G|^{-1}$ exists in $k$,
then for any irreducible character
$\chi$ of $H$, $\chi(\nu_H)$ is identical to Bantay's formula
(\ref{eq: bantayfsind}). We begin with the definition of twisted
quantum doubles of finite groups. \\

Let $G$ be a finite group and $\omega: G \times G \times G \map k^\times$
be a normalized 3-cocycle; that is, a function such that
$\omega(x,y,z)=1$ whenever one of $x,y$ or $z$ is equal to the
identity element 1 of $G$ and which satisfies the functional equation
\begin{equation}\label{eq:3cocycle}
\omega(g,x,y)\omega(g, xy,
z)\omega(x,y,z)=\omega(gx,y,z)\omega(g,x,yz)\quad \mbox{for
any } g, x, y, z \in G\,.
\end{equation}
For any $g \in G$, define the functions $\theta_g, \gamma_g :
G\times G \rightarrow k^\times$ as follows:
\begin{eqnarray}
\theta_g(x,y) &=&\frac{\omega(g,x,y)\omega(x,y,(x y)^{-1}g x
y)}{\omega(x,x^{-1}g x,y)}\,,\label{eq0.001}\\
\gamma_g(x,y) &=& \frac{\omega(x,y,g)\omega(g, g^{-1}x g,
g^{-1}yg)}{\omega(x,g,
g^{-1}y g)}\,.\label{eq0.002}
\end{eqnarray}
Let $\{e(g)|g \in G\}$ be the dual basis of the canonical basis of
$k [G]$. The {\em twisted quantum double} $D^\w(G)$ of $G$ with
respect to $\omega$ is the quasi-Hopf algebra with
underlying vector space $k[G]' \otimes k[G]$. The
multiplication, comultiplication and associator are given,
respectively, by
\begin{equation}\label{multiplication}
(e(g) \otimes x)(e(h) \otimes y) =\theta_g(x,y)
  \delta_{g, xhx^{-1}}e(g)\otimes x y\,,
\end{equation}
\begin{equation}\label{comultiplication}
\Delta(e(g)\otimes x)  = \sum_{hk=g} \gamma_x(h,k) e(h)\otimes x
\otimes e(k) \otimes x\,,
\end{equation}
\begin{equation}\label{eq0.01}
\Phi = \sum_{g,h,k \in G} \omega(g,h,k)^{-1} e(g) \otimes 1 \otimes
e(h) \otimes 1 \otimes e(k) \otimes 1\,.
\end{equation}
The counit and antipode are given by
\begin{equation}\label{eq:counit}
  \varepsilon(e(g)\otimes x) = \delta_{g,1}
\end{equation}
and

\begin{equation}\label{eq:antipode}
   S(e(g)\otimes x) =
\theta_{g^{-1}}(x,x^{-1})^{-1}\gamma_x(g,g^{-1})^{-1}e(x^{-1}g^{-1}x)\otimes
x^{-1}\,,
\end{equation}
where $\delta_{g,1}$ is the Kronecker delta. The corresponding
elements $\alpha$ and $\beta$ are $1_{D^{\w}(G)}$ and $\sum\limits_{g \in
G} \w(g,g^{-1},g)e(g)\otimes 1$ respectively (cf. \cite{DPR90}).
Verification of the
detail involves the following identities, which result from the
3-cocycle identity for $\w$:
\begin{equation}\label{eq:0.1}
\theta_z(a,b)\theta_z(a b,c) = \theta_{a^{-1}z a}(b,c)\theta_z(a,b
c)\,,
\end{equation}
\begin{equation}\label{eq:0.2}
\theta_y(a,b)\theta_z(a,b)\g_a(y,z)\g_b(a^{-1}y a,a^{-1}z
a)=\theta_{yz}(a,b)\g_{ab}(y,z)\,,
\end{equation}
\begin{equation}\label{eq:0.3}
\g_z(a,b)\g_z(a b,c)\w(z^{-1}a z,z^{-1}b z,z^{-1}c
z)=\g_z(b,c)\g_z(a,b c)\w(a,b,c)\,,
\end{equation}
for all $a,b,c,y,z \in G$.

\begin{remark}\label{r0.1}
{\rm The algebra $D^\w(G)$ is a semi-simple (cf. \cite{DPR90}).
If $\omega=1$, then the twisted quantum double $D^{\omega}(G)$
identical to the Drinfeld double of the group algebra $k [G]$.
However, $D^{\omega}(G)$ is not a Hopf algebra in general. Moreover,
even if $\omega, \omega'$ differ by a coboundary,  $D^\w(G)$ and
$D^{\w'}(G)$ are not isomorphic as quasi-bialgebras. Nevertheless,
they are {\em gauge equivalent}. In addition, if $G$ is abelian, $D^\w(G)$
also admits a Hopf algebra structure with the same underlying
$\Delta$, $\e$
and $S$ (cf. \cite{MN01}).
}
\end{remark}

Let
\begin{equation}\label{integral}
\Lambda = \frac{1}{|G|}\sum_{x \in G} e(1) \otimes x \in D^\omega(G)\,.
\end{equation}
It is straightforward
to show that $\Lambda$ is a left integral of $D^\omega(G)$. Moreover,
$$
\e(\Lambda)=1  \,.
$$
After \cite{Pana98} and \cite{HNQA}, this gives another proof of the
semi-simplicity of $D^\omega(G)$.
Note that
$$
\Delta(\Lambda) = \sum \Lambda_1 \otimes \Lambda_2  = \frac{1}{|G|}
\sum_{g,x \in G}
\gamma_x(g, g^{-1})e(g) \otimes x \otimes e(g^{-1}) \otimes x \,.
$$
Since $\beta \alpha = \beta$ is invertible, it follows from Corollary
\ref{cor3.5} that
\begin{eqnarray*}
\nu_{D^\w(G)} &=&
\frac{1}{|G|}\left(\sum_{g \in G} \omega(g, g^{-1}, g)^{-1}(e(g)
\otimes 1)\right) \left(\sum_{g,x \in G}
\gamma_x(g, g^{-1})(e(g) \otimes x) (e(g^{-1}) \otimes x)\right)\\
&=& \frac{1}{|G|} \left(\sum_{g \in G} \omega(g^{-1}, g, g^{-1})(e(g)
\otimes 1)\right)
\left(\sum_{x^{-1}gx=g^{-1}}\gamma_x(g, g^{-1})\theta_g(x,x)(e(g)
\otimes x^2)\right)
\\
&=& |G|^{-1} \sum_{x^{-1}gx=g^{-1}}
\omega(g^{-1}, g, g^{-1})\gamma_x(g, g^{-1})\theta_g(x,x)(e(g) \otimes x^2).
\end{eqnarray*}
Here we have used the equality
$$
\w(g, g^{-1}, g)^{-1}= \w(g^{-1}, g, g^{-1})
$$
which is readily derived from equation (\ref{eq:3cocycle}). Thus for any
irreducible character $\chi$ of $D^\w(G)$, the Frobenius-Schur indicator of
$\chi$ is
$$
\chi(\nu_{D^\w(G)})=|G|^{-1} \sum_{x^{-1}gx=g^{-1}}
\omega(g^{-1}, g, g^{-1})\gamma_x(g, g^{-1})\theta_g(x,x)\chi(e(g) \otimes x^2)
$$
as given by Bantay.

\section{Trace Elements and Antipodes of Semi-simple Quasi-Hopf Algebras}
It is proved by Larson and Radford \cite{LaRa87} \cite{LaRa88}
that if char $k = 0$, the antipode of a semi-simple Hopf algebra over $k$ is an
involution. However, the antipode of a semi-simple quasi-Hopf
algebra $H$ could be of any order. Nevertheless,  we prove an
analog of the Larson-Radford theorem for a split semi-simple quasi-Hopf algebras
$H$ over any field $k$: there exists a unit $u \in H$ such that the
antipode of $H_u$ is an involution. To this end we introduce the
\emph{trace element} $g$ of a semi-simple quasi-Hopf algebra. This
element will play a role throughout the remaining Sections of the
paper. \\

Let $H=(H, \Delta,\e, \Phi, \a, \b, S)$ be a finite-dimensional semi-simple
quasi-Hopf algebra
over $k$ and $\Lambda$ the normalized two-sided integral of $H$.
By \cite{HNQA}, there exists a functional $\lambda \in H'$, called
the normalized
left cointegral of $H$, given by the formula
\begin{equation}\label{eq:cointegral}
\lambda(x)=\sum_i b^i(x S^2(b_i)S(\b)\a)
\end{equation}
for all $x \in H$, where $\{b_i\}$ is a basis of $H$ and $\{b^i\}$ is its dual
basis (see \cite{HNQA} for the details of cointegral). The normalized
left cointegral $\lambda$ admits the following properties :
\begin{enumerate}
\item[(i)] $\lambda(\Lambda)=1$.
\item[(ii)] $\lambda(ab)$ ($a, b \in H$) defines a non-degenerate bilinear
form on $H$.
\item[(iii)] For all $a, b \in H$,
\begin{equation}\label{eq:nakayama}
\lambda(ab)=\lambda(bS^2(a)) \,.
\end{equation}
\end{enumerate}

Let $\chi_{reg}$ denote the character of the left regular representation of $H$.
The bilinear form on $H$ defined by $\<a,b\>_{reg}:=\chi_{reg}(ab)$ is
then  symmetric
and non-degenerate. By the non-degeneracy of $\lambda$,
there exists a unique element $g$ of $H$ such that
\begin{equation}\label{eq:trace_element}
 \chi_{reg}(x) = \lambda(xg)
\end{equation}
for all $x \in H$. We call $g$ the {\em trace element}.

\begin{example}
{\rm
If char $k =0$,  and $H$ is a finite-dimensional semi-simple Hopf algebra over $k$, then $S^2=id_H$.
By (\ref{eq:cointegral}),
$$
\lambda(x)=\sum_i b^i(xb_i)=\chi_{reg}(x)\,.
$$
Thus, the trace element of $H$ is 1.
}
\end{example}

\begin{lem}\label{l:inner}
Let $H=(H, \Delta,\e, \Phi, \a, \b, S)$ be a finite-dimensional semi-simple
quasi-Hopf algebra.
Then the trace element $g$ of $H$ is invertible and
$$
S^2(a) =g^{-1}ag
$$
for all $a \in H$. Moreover, $gS(g)$ is in the center of $H$ and $gS(g)=S(g)g$.
\end{lem}
\pf By (\ref{eq:trace_element}), the left annihilator of $g$ in $H$
is a subset of $\ker \chi_{reg}$. Since $H$ is semi-simple, $\ker \chi_{reg}$
does not contain any non-trivial left ideals of $H$.
Therefore, the left annihilator of $g$ is trivial. Since the left regular
representation of $H$ is faithful and finite-dimensional, $g$ is invertible.
Thus, we have
$$
\lambda(ab)=\lambda(abg^{-1}g)= \chi_{reg}(abg^{-1})=\chi_{reg}(bg^{-1}a)=\lambda(bg^{-1}ag)
$$
for all $x, y \in H$. By the non-degeneracy of $\lambda$ and
(\ref{eq:nakayama}), we obtain
$$
S^2(a)=g^{-1}ag
$$
for all $a \in H$. In particular,
$$
S(g^{-1}ag)=S^3(a)=g^{-1}S(a)g\,.
$$
Therefore,
\begin{equation} \label{eq:g}
gS(g)S(a)=S(a)gS(g)
\end{equation}
for all $a \in H$ and hence $gS(g)$ is in the center of $H$. Taking $a=g^{-1}$
in (\ref{eq:g}), the result in the last statement follows. \qed\\

\begin{lem}\label{l:u}
Let $A$ be a finite-dimensional split semi-simple algebra over $k$
and $S$ an algebra anti-automorphism on $A$ such that $S^2$ is inner. Then
there exists a unit $u \in A$ such that $S_u^2 =id_A$ where
$$
S_u(x) = uS(x)u^{-1}
$$
for all $x \in A$.
\end{lem}
\pf Without loss of generality, we can assume that $A$ is a direct
sum of full matrix rings over $k$, say $A=\oplus_{i=1}^d
M_{n_i}(k)$.  Let $\iota_i $ denote the natural embedding from
$M_{n_i}(k)$ into $A$, $p_i$  the natural surjection from $A$ onto
$M_{n_i}(k)$, and $A_i$ the image of $\iota_i$. Then, $A_1, \cdots
, A_d$ is the complete set of minimal ideals of $A$. Since $S$ is
an algebra anti-automorphism, there exists a permutation $\sigma$
on $\{1,\dots, d\}$ such that $S(A_i) = A_{\sigma(i)}$ for all
$i=1, \dots, d$. As $S^2$ is inner, $S^2(A_i)=A_i$ for all $i$ and
so $\sigma^2=id$.\\

Since $S(A_i)= A_{\sigma(i)}$, $M_{n_i}(k) =
M_{n_{\sigma(i)}}(k)$. Moreover, $p_j \circ S \circ \iota_i=0$ for
$j \ne \sigma(i)$ and $p_{\sigma(i)} \circ S \circ \iota_i$ is  an
algebra anti-automorphism on $M_{n_i}(k)$. By the Skolem-Noether
theorem, there exists an invertible matrix $u_{i} \in
M_{n_{\sigma(i)}}(k)$ such that $p_{\sigma(i)} \circ  S\circ \iota
(x) =u_i^{-1}x^t u_i$ for any $x \in M_{n_i}(k)$ where $x^t$ is
the transpose of $x$.\\

Let $u=\sum_{i=1}^d \iota_{\sigma(i)}(u_i)$. Since $u_i$ is
invertible in $M_{n_{\sigma(i)}}(k)$ for all $i$, $u$ is
invertible in $A$. Since $S(A_i) = A_{\sigma(i)}$ is an ideal of
$A$, $S_u(A_i) = A_{\sigma(i)}$. Then for any $x \in M_{n_i}(k)$,
$$
p_{\sigma(i)}(S_u(\iota_i(x))) = u_i(p_{\sigma(i)}\circ S \circ
\iota_i(x))u_i^{-1}=x^t\,.
$$
Thus,
$$
\iota_{\sigma(i)}(x^t) = \iota_{\sigma(i)}\circ
p_{\sigma(i)}(S_u(\iota_i(x))) = S_u(\iota_i(x))\,,
$$
and hence
$$
S^2_u(\iota_i(x))=S_u(\iota_{\sigma(i)}(x^t))=\iota_{\sigma^2(i)}((x^t)^t)
=\iota_i(x)
$$
as $\sigma^2=id$. Therefore, $S^2_u(a)=a$ for all $a \in A_i$,
$i=1, \dots, d$. Since $A=A_1 \oplus \cdots \oplus A_d$,
$S^2_u=id_A$. \qed

\begin{thm}\label{t:6.4}
Let $H=(H, \Delta,\e, \Phi, \a, \b, S)$ be a finite-dimensional split semi-simple quasi-Hopf
algebra over $k$. Then there exists an invertible element $u$ of $H$
such that the antipode of $H_u$ is an involution.
\end{thm}
\pf It follows from Lemma \ref{l:inner} or \cite[Proposition 5.6]{HNQA}
 that $S^2$ is inner. By Lemma \ref{l:u}, the result follows. \qed\\

\begin{remark} \label{r:equiv_dual}
{\rm
Suppose $u$ is an invertible element of $H$ and
$M$  a finite-dimensional left $H$-module. Let  ${^+\!M}$,  $^*\!M$ denote  the left dual of
$M$ in $H_u$-$\mod_{fin}$ and $H$-$\mod_{fin}$ respectively. Then, ${^+\!M}$ and $^*\!M$ are
isomorphic left $H$-modules under the map $\phi_u : {^+\!M} \map {^*\!M}$ defined
by
$$
\phi_u(f)(x)=f(ux)
$$
for all $x \in M$ and $f \in M'$. In particular, $M \cong {^*\!M}$ if, and
only if, $M \cong {^+\!M}$ as left $H$-modules (cf. \cite[p1425]{Drin90}).
}\qed
\end{remark}
\section{Pivotal Category Structure of {$H\mbox{-}\mod_{fin}$}}
We begin (Theorem \ref{t:etingof}) with Etingof's observation that the trace element
 $g$  of a
finite-dimensional semi-simple quasi-Hopf algebra $H=(H, \Delta,\e, \Phi, \a, \b, S)$
over an algebraically closed field of characteristic zero
defines an isomorphism of tensor functors
$$
 j: Id \map \,{^{**}?}\,.
$$
Moreover, we prove that $S(g)=g^{-1}$, a fact that we will need in
Section 8. A direct result of this is that
$H$-$\mod_{fin}$ is a pivotal category in the sense of Joyal and Street
(cf. \cite{FY}). For the remainder of this paper we will assume
that $k$ is an algebraically closed field of characteristic zero.\\

For simplicity, we write $\mathcal{C}$ for the semi-simple rigid
tensor category $H$-$\mod_{fin}$ in this section. Obviously,
$\mathcal{C}$ is a fusion category over $k$ (cf. \cite{ENO}).
Recall from \cite{BaKi} that if $V \in \mathcal{C}$ and $f: V
\map\, {^*V}$ then the categorical trace of $f$ is the scalar
$\tr_V(f)$ defined by
\begin{equation}\label{e:dim}
\ev_{^*\!V} \circ (f \otimes id) \circ \coev_V\,.
\end{equation}
where $\ev_V: {^*V} \otimes V \map k$ and
$\coev_V: k\map\, V \otimes  {^*V}$ are evaluation and coevaluation maps.

Following \cite{Mu}, for any simple object $V$ in $\mathcal{C}$ and an
isomorphism $f: V \map\, {^*V}$, we define
\begin{equation}\label{e:squared_norm_1}
|V|^2=\tr_V(f)\,\tr_{^*\!V}({^*\!(f^{-1})})\,.
\end{equation}
Clearly, $|V|^2$ is independent of the choice of $f$.\\

By \cite{ENO}, there exists an isomorphism of tensor functors
$$
j: Id \map \,{^{**}?}
$$
such that for any simple object $V$ of $\mathcal{C}$,
\begin{equation}\label{e:dimensions}
\tr_V(j) = \FPdim(V)=\dim(V)\,
\end{equation}
where $\FPdim(V)$ is the {\em Frobenius-Perron dimension} of $V$. Moreover,
\begin{equation}\label{e:squared_norm_2}
|V|^2 = \dim(V)^2\,.
\end{equation}

Let $a$  be the unique invertible element of $H$  such that
\begin{equation}\label{e:a}
j_H(1)(f)=f(a)
\end{equation}
for all $f \in {^*\!H}$. By the naturality of $j$, one can show that
\begin{equation}\label{e:inner_a}
 S^2(x)=axa^{-1} \quad \mbox{ for all } x \in H\,,
\end{equation}
 and for any $V \in \mathcal{C}$, $j:V \map\, {^{**}V}$
is given by
\begin{equation}\label{e:j}
j_V(x)(f)= f(ax)
\end{equation}
for all $x \in V$ and $f \in {^*V}$. Thus, by (\ref{e:dim}) and (\ref{e:squared_norm_1}),
for any simple objective $V$ in $\mathcal{C}$ with character $\chi$,
\begin{equation}\label{e:dimV}
\dim(V)= \chi(a \b S(\a))\,,
\end{equation}
and
\begin{equation}\label{e:squared norm}
|V|^2= \chi(a \b S(\a))\chi(a^{-1} S(\b) \a)\,.
\end{equation}
Hence, by (\ref{e:dimV}) and (\ref{e:squared_norm_2}), we also have
\begin{equation}\label{e:dimV_2}
\dim(V)=\chi(a^{-1} S(\b) \a)
\end{equation}
In fact, $a^{-1}$ is the trace element of $H$.\\
\begin{thm}\label{t:etingof}
Let $H=(H, \Delta,\e, \Phi, \a, \b, S)$ be a finite-dimensional semi-simple
quasi-Hopf algebra over $k$ and $g$ the trace element of $H$. Then the
natural isomorphism $j_V: V \map\, {^{**}\!V}$ for any $V$ in $H$-$\mod_{fin}$,
given by
$$
j_V(x)(f)=f(g^{-1}x)
$$
for all $x \in V$ and $f \in {^{*}\!V}$, defines an isomorphism of the tensor
functors $Id$ and ${^{**}?}$ such that
$$
\dim(V) = \chi(g^{-1} \b S(\a))
$$
for any simple $H$-module  $V$ with character $\chi$.
\end{thm}
\pf By the preceding discussion, it suffices to show that the
element $a$ defined in (\ref{e:a}) is identical to $g^{-1}$.
By Lemma \ref{l:inner} and (\ref{e:inner_a}), $ag$ is in the center of $H$.
Therefore, it is enough to show that for any simple $H$-module $V$
with character $\chi$,
$$
\chi(a \b S(\a) ) =\chi( g^{-1} \b S(\a) )\,.
$$
Let $e_V$ be the central idempotent of $H$ such that
$$
\chi(x)\dim(V)=\chi_{reg}(e_V x)
$$
for all $x \in H$. Thus, we obtain
\begin{eqnarray*}
\chi(g^{-1}\b S(\a))\dim(V) & = & \chi_{reg}(e_V g^{-1}\b S(\a)) \\
&=& \chi_{reg}(e_V \b S(\a) g^{-1}) \\
&=& \lambda(e_V \b S(\a))
\end{eqnarray*}
where $\lambda$ is the normalized left cointegral of $H$. Let $\{b^i\}$ be the dual basis of
the basis $\{b_i\}$ of $H$. Then, we have
\begin{eqnarray*}
\chi(g^{-1}\b S(\a))\dim(V) & = & \sum_i b^i(e_V \b S(\a)S^2(b_i) S(\b)\a) \\
&=& \chi_{reg}(e_V \b S(\a) a b_i a^{-1} S(\b)\a) \\ \\
&=& \chi(\b S(\a) a)\chi(a^{-1} S(\b)\a) \\ \\
&=& |V|^2 = \dim(V)^2 \quad \mbox{by (\ref{e:squared norm}) and
(\ref{e:squared_norm_2})}\,.
\end{eqnarray*}
Therefore, by (\ref{e:dimV}), we obtain
$$
\chi(g^{-1}\b S(\a))=\dim(V)=\chi(a \b S(\a))\,.
$$
\mbox{}\hfill\qed

\begin{thm}\label{t:7.2}
Let $H=(H, \Delta,\e, \Phi, \a, \b, S)$ be a finite-dimensional
semi-simple quasi-Hopf algebra over $k$ and $g$ the trace element
of $H$. Then $S(g)=g^{-1}$ and hence
\begin{equation}\label{e:pivotal}
{^*(j_V)}\circ j_{^*V} =id_{^*V}
\end{equation}
for any $V \in H$-$\mod_{fin}$
\end{thm}
\pf Since $gS(g)$ is central, $gS(g)$ acts on any simple
$H$-module $V$ as multiplication by a scalar $c_V\in k$. In order
to show that $S(g)=g^{-1}$, it suffices to prove that
$$
c_V=1
$$
for any simple $H$-module $V$.\\

Let $V$ be a simple $H$-module with character $\chi$. Then the
character of ${^*V}$ is ${^*\!\chi}$ given by
$$
{^*\!\chi}=\chi\circ S\,.
$$
By Theorem \ref{t:etingof},  (\ref{e:dimV}) and (\ref{e:dimV_2}),
we have
$$
\begin{aligned}
 \dim({^*V})& = {^*\!\chi}(gS(\b)\a)=
 \chi(S(\a)S^2(\b)S(g))\\
 &= \chi(S(\a) g^{-1} \b gS(g)) = c_V \, \chi(S(\a) g^{-1}
 \b)\\
  & = c_V\, \dim(V)
\end{aligned}
$$
Therefore, $c_V=1$. Equation (\ref{e:pivotal}) follows easily from $S(g)=g^{-1}$.
\qed\\

Theorem \ref{t:etingof} and (\ref{e:pivotal}) implies that $H$-$\mod_{fin}$ is indeed
a pivotal category defined by Joyal-Street (cf. \cite{FY}). Nikshych
also pointed out that (\ref{e:pivotal}) can be proved using \emph{weak Hopf algebras}.

\section{Frobenius-Schur Indicators via Bilinear Forms with Adjoint $S$}
Let $H=(H, \Delta,\e, \Phi, \a, \b, S)$ be a finite-dimensional semi-simple quasi-Hopf algebra
over $k$, and $g$ the trace element of $H$. In this section,
we will prove
that for any simple left $H$-module $M$ with character $\chi$, the Frobenius-Schur
indicator $\chi(\nu_H)$ of $\chi$ can only be 0, 1 or -1.
It is non-zero if, and only if $M \cong {^*\!M}$. Moreover in this case, $M$ admits a non-degenerate bilinear form
 $\<\cdot, \cdot\>$
such that $\<h u, v\> = \< u, S(h)v\>$ for all $h \in H$, $u, v \in M$,
and
$$
\<u, v\>  =\chi(\nu_H) \<v, g^{-1}u\>\,.
$$

\begin{defn}
Let $H=(H, \Delta,\e, \Phi, \a, \b, S)$ be a quasi-Hopf algebra
over $k$,  $M$ be a left $H$-module and
$\<\cdot, \cdot\>$ a bilinear form on $M$.
\begin{enumerate}
\item[\rm (i)] The form is said to be $H$-invariant if
$$
\sum \<h_1 u , h_2 v\> = \e(h) \<u,v\>
$$
for all $h \in H$ and $u,v \in V$  where $\sum h_1 \otimes h_2 = \Delta(h)$.
\item[\rm (ii)] The antipode $S$ is said to be the adjoint of the form if
$$
 \<h u ,  v\> =  \<u, S(h) v\>
$$
for all $h \in H$ and $u,v \in V$.
\end{enumerate}
\end{defn}

\begin{lem}\label{l:schur}
Let $H=(H, \Delta,\e, \Phi, \a, \b, S)$ be a quasi-Hopf algebra over
$k$ and  $M$  a simple left $H$-module. If $\<\cdot, \cdot\>_1$
and $\<\cdot, \cdot\>_2$ are non-degenerate bilinear forms on $M$
with the same adjoint $S$, then there exists a non-zero element $c \in
k$ such that
$$
\<u, v\>_1 = c \<u, v\>_2
$$
for all $u, v \in M$.
\end{lem}
\pf Define  $J_i : M \map {^*\!M}$ ( $i=1, 2$) by
$$
J_i(u)(v) = \< u, v\>_i
$$
for $u, v \in M$. Since $\<\cdot, \cdot\>_1$ and $\<\cdot,
\cdot\>_2$ are non-degenerate bilinear forms on $M$ with the
adjoint $S$, $J_1, J_2$ are isomorphisms of $H$-modules. In
particular, $M$ and ${^*\!M}$ are isomorphic simple $H$-modules.
By Schur's lemma, $J_1=c J_2$ for some non-zero element $c \in k$
and so the result follows. \qed\\

\begin{lem}\label{l:dualbasis}
Let $H=(H, \Delta,\e, \Phi, \a, \b, S)$ be a finite-dimensional semi-simple quasi-Hopf algebra, $\Lambda$ the
normalized two-sided integral of $H$ and $g$ the trace element of $H$.
Suppose that
$$
q_R \Delta(\Lambda) p_R = \sum_{i=1}^n x_i \otimes y_i
$$
where $\{x_i\}$ is basis of $H$. Then $\{S(x_i) g^{-1}, y_i\}$ is a pair
dual bases with
respect to $\<\cdot,\cdot\>_{reg}$.
\end{lem}
\pf Following \cite{HNQA}, we define the elements $U, V \in H
\otimes H$ by
\begin{eqnarray}
U &= &F_H^{-1}(S \otimes S)(q_R^{21})\,,\\
V& = &(S^{-1} \otimes S^{-1})(F_H^{21}p_R^{21})
\end{eqnarray}
where $F_H$, $q_R, p_R \in H \otimes H$ are defined in (\ref{FH})
and (\ref{R}). By \cite[(7.3) and (7.4)]{HNQA},
\begin{eqnarray*}
q_R \Delta(\Lambda) p_R &=&
(q_L^2 \otimes 1)V\,\Delta(S^{-1}(q_L^1))\, \Delta(\Lambda)
\Delta(S(p_L^1))\,U(p_L^2 \otimes 1)\\
&=& (q_L^2 \e(S^{-1}(q_L^1)) \otimes 1)V\,\Delta(\Lambda)\,U (\e(S(p_L^1))p_L^2 \otimes 1)\\
\end{eqnarray*}
By \cite[Remark 7]{Drin90}, $\e \circ S = \e = \e \circ S^{-1}$. Therefore,
$$
q_L^2 \e(S^{-1}(q_L^1)) = \e(\a)1_H, \quad \mbox{and}\quad \e(S(p_L^1))p_L^2=\e( \b)1_H\,.
$$
It follows from (\ref{2.9}) that $\e(\a\b)=1$ and so
$$
q_R \Delta(\Lambda) p_R = \e(\a)\e(\b) V \Delta(\Lambda) U =V \Delta(\Lambda) U \,.
$$
Let $\lambda$ be the normalized left cointegral of $H$.
By \cite[Proposition 5.5]{HNQA},
$$
\sum_i S(x_i)\lambda(y_i a) = a
$$
for all $a \in H$. In particular,
$$
a = (ag)g^{-1} = \sum_i S(x_i)g^{-1} \lambda(y_i ag) =
\sum_i S(x_i)g^{-1} \chi_{reg}(y_i a)\,.
$$
Since $\{S(x_i)g^{-1}\}$ is also a basis of $H$,
$\chi_{reg}(y_i S(x_j)g^{-1}) = \delta_{ij}$ and so $\{S(x_i) g^{-1}, y_i\}$ is a pair
dual bases of $H$ with respect to $\<\cdot, \cdot\>_{reg}$. \qed\\

\begin{lem}\label{l:sum_basis}
Let $A$ be a finite-dimensional semi-simple algebra over $k$ and
$\{a_i, b_i\}$ a pair dual bases with respect to the form
$\<\cdot, \cdot\>_{reg}$. Then
$$
   \sum_i a_i b_i =1_A\,.
$$
\end{lem}
\pf Without loss of generality, we may assume that $A=\oplus_{i=1}^d
M_{n_i}(k)$. Then $\chi_{reg}(x)= \sum_{i=1}^d n_i tr_i(x)$ where
$tr_i(x)$ is the trace of the $i$th component matrix of $x$. Let
$\{e^i_{lm}\}$ be the set of matrix units for the $i$th summand
$M_{n_i}(k)$ of $A$. Following  \cite{LM00},
$\{n_i^{-1}e^i_{lm}, e^i_{ml}\}$ is  a pair of dual basis with respect to
$\<\cdot, \cdot\>_{reg}$. Thus,
$$
\sum_{i, l, m} n_i^{-1}e^i_{lm}e^i_{ml}= \sum_{i, l}e^i_{ll} =
1_A\,.
$$
It follows from \cite[Lemma 2.6]{LM00} that
$$
\sum_i a_i b_i = \sum_{i, l,m} n_i^{-1}e^i_{lm}e^i_{ml} =1_A \,.
$$
\hfill\qed\\
\begin{cor}
Let $H=(H, \Delta,\e, \Phi, \a, \b, S)$ be a finite-dimensional
semi-simple  quasi-Hopf algebra over $k$. Then trace element $g$ of $H$ is
given by
$$
g = m\tau(S \otimes id)(q_R\Delta(\Lambda) p_R)
$$
where $\Lambda$ is the normalized integral of $H$, $m$ is  multiplication and $\tau$
 the usual flip map.
\end{cor}
\pf Let
$$
q_R \Delta(\Lambda) p_R = \sum_i x_i \otimes y_i\,.
$$
By Lemma \ref{l:dualbasis} and Lemma \ref{l:sum_basis}, we have
$$
\sum_i y_iS(x_i)g^{-1} =1
$$
and so the result follows. \qed\\

Let $\{a_i, b_i\}$ be  dual bases of the semi-simple quasi-Hopf algebra $H$
with respect to the form $\< \cdot, \cdot\>_{reg}$ discussed in Lemma
\ref{l:dualbasis}.
For any $k$-involution $\mathcal{I}$ on $H$ and
for any character $\chi$ of $H$, we define
$$
\mu_2(\chi, \mathcal{I})= \chi(\sum_i \mathcal{I}(a_i)b_i)\,.
$$
\begin{remark} \label{r:mu2}
{\rm
Since  $\sum_i a_i b_i = 1_H$ by Lemma \ref{l:sum_basis},
the $\mu_2$ defined in \cite[Therorem 2.7]{LM00} with respect to the
$k$-involution $\mathcal{I}$ is given by
$$
\frac{\chi(1_H)}{\chi(\sum_i a_i b_i)}
\chi(\sum_i \mathcal{I}(a_i)b_i) = \chi(\sum_i \mathcal{I}(a_i)b_i)
$$
which coincides with $\mu_2(\chi, \mathcal{I})$.
}
\end{remark}

\begin{lem}\label{l:mu2}
Let $H=(H, \Delta,\e, \Phi, \a, \b, S)$ be a finite-dimensional
semi-simple  quasi-Hopf algebra over $k$, $g$ the trace element of
$H$, and $M$ an irreducible $H$-module with character $\chi$. Then for any
unit $u \in H$ such that $S_u$ is an involution,
$$
\mu_2(\chi, S_u)= c\,\chi(\nu_H)
$$
where $c$ is the non-zero scalar given by
$$
c=\frac{\chi( u S(u^{-1})g^{-1})}{\dim M}\,.
$$
\end{lem}
\pf If $u$ is a unit of $H$ such that $S_u$ is an involution, then
for any $x \in H$,
$$
x= S_u^2(x)=uS(u^{-1}) S^2(x) S(u) u^{-1}
$$
or equivalently
$$
S^2(x)=S(u)u^{-1} x u S(u^{-1})\,.
$$
By Lemma \ref{l:inner}, $ uS(u^{-1})g^{-1}$ is in the center of $H$.
Thus, $uS(u^{-1}) g^{-1}$ acts
on $M$ as multiplication by the non-zero scalar
$$
c=\frac{\chi(uS(u^{-1})g^{-1}))}{\dim M}\,.
$$
Suppose that
$$
q_R \Delta(\Lambda) p_R = \sum_i x_i \otimes y_i
$$
as in Lemma \ref{l:dualbasis} where $\Lambda$ is the normalized two-sided integral
of $H$. Then we have
\begin{eqnarray*}
\mu_2(\chi, S_u) & = &  \chi(\sum_i S_u(S(x_i)g^{-1}) y_i) \\
            &=& \chi(\sum_i u S(g^{-1})S^2(x_i) u^{-1} y_i) \\
            &=& \chi(\sum_i u S(g^{-1}) g^{-1} x_i g u^{-1} y_i) \\
            &=&\chi(\sum_i u S(g^{-1})g^{-1}S(u^{-1})S(g) x_i y_i) \quad\mbox{(by Lemma \ref{pq_1})}\\
            &=& \chi(\sum_i u S(u^{-1})g^{-1}x_i y_i) \quad
            \mbox{(by Lemma \ref{l:inner})}\\\\
            & =& c\, \chi(\nu_H)\,.
\end{eqnarray*}
\begin{thm} \label{t:main2}
Let $H=(H, \Delta,\e, \Phi, \a, \b, S)$ be a finite-dimensional
semi-simple  quasi-Hopf algebra over $k$, $g$ the trace element of
$H$, and $M$ a simple
$H$-module with character $\chi$. Then the Frobenius-Schur indicator
$\chi(\nu_H)$ of $\chi$ satisfies the following properties:
\begin{enumerate}
\item[\rm (i)] $\chi(\nu_H)  \ne 0$ if, and only if, $M \cong {^*\!M}$ as
left $H$-modules.
\item[\rm (ii)]  For any non-zero $\kappa \in k$,
$\chi(\nu_H) =\kappa$ if, and only if, $M$ admits a
non-degenerate bilinear form $\< \cdot, \cdot\>$ with
the adjoint $S$ such that
$$
\< x, y \> = \kappa\<y, g^{-1}x \>
$$
for all $x, y \in M$.
\item[\rm (iii)] The values of $\chi(\nu_H)$ can only be $0$, $1$ or
$-1$.
\end{enumerate}
Moreover,
$$
 \Tr(S) =  \sum_{\chi \in Irr(H)} \chi(\nu_H)\chi(g^{-1}).
$$
\end{thm}
\pf By Theorem \ref{t:6.4}, there exists an unit $u \in H$ such that
$S_u$ is an involution. As in the proof of Lemma \ref{l:mu2},
$uS(u^{-1}) g^{-1}$  is a central unit of $H$. Thus, $uS(u^{-1}) g^{-1}$ acts
on $M$ as multiplication by the non-zero scalar
$$
c=\frac{\chi(uS(u^{-1})g^{-1})}{\dim M}\,.
$$
Also, by \cite[Theorem 2.7]{LM00} and Remark \ref{r:mu2}, the element
$\mu_2(\chi, S_u) \ne 0$ if, and only if
$M \cong {^+\!M}$ as left $H$-modules where ${^+\!M}$ is the left $H$-module
with underlying space $M'$ and the $H$-action given by
$$
(hf)(x) = f(S_u(h)x)
$$
for all $f \in M'$ and $h \in H$. Actually, ${^+\!M}$ is the left
dual of $M$ in
$H_u$-$\mod_{fin}$. It follows from  Remark \ref{r:equiv_dual} that $\mu_2(\chi, S_u) \ne 0$
if, and only if $M \cong {^*\!M}$ as left $H$-modules.
Hence, by Lemma \ref{l:mu2}, statement (i) follows.\\

If $\chi(\nu_H) \ne 0$, then $\mu_2(\chi, S_u) \ne 0$ by Lemma \ref{l:mu2}.
By Remark \ref{r:mu2} and
\cite[Theorem 2.7(ii)]{LM00}, $M$ admits a non-degenerate
bilinear form
$(\cdot, \cdot)$ with adjoint $S_u$ such that
$$
(x, y) = \mu_2(\chi, S_u) (y,x)
$$
for any $x, y  \in M$. Define
$$
\<x, y\>= (x, uy)
$$
for any $x, y \in M$. One can easily see that $\<\cdot, \cdot\>$ is a
non-degenerate bilinear form on $M$ with adjoint $S$. Moreover, for any
$x, y\in M$,
$$
\<x, y\> = (x, uy) = \mu_2(\chi, S_u) (uy,x)=  \mu_2(\chi, S_u) (y,S_u(u) x)\,.
$$
Thus, by Lemma \ref{l:mu2}, we obtain
$$
 \<x, y\> =c\, \chi(\nu_H) (y, S_u(u)x)
=  \chi(\nu_H) \<y, S(u)u^{-1}cx\> = \chi(\nu_H) \<y, g^{-1}x\>\,.
$$

Conversely, suppose $M$ admits  a non-degenerate bilinear form
$\<\cdot, \cdot\>$ with adjoint $S$ and that there  exists a non-zero
element $\kappa$ of  $k$ such that
$$
\< x, y \> = \kappa\<y, g^{-1} x \>
$$
for all $x, y \in M$. Then the map $J: M \map {^*\!M}$,
defined by
$$
J(x)(y)=\<x, y\>, \quad x, y \in M\,,
$$
 is an isomorphism of left $H$-modules. Thus, by (i),
$\chi(\nu_H)\ne 0$. Hence, by above arguments, $M$ admits a
non-degenerate bilinear form
$\<\cdot, \cdot\>_0$ with adjoint $S$ such that
$$
\< x, y \>_0 = \chi(\nu_H) \<y, g^{-1} x \>_0
$$
for all $x, y \in M$. By Lemma \ref{l:schur}, $\<\cdot, \cdot\>$ is a
non-zero scalar multiple of $\<\cdot, \cdot\>_0$. Therefore,
$$
\kappa=\chi(\nu_H)
$$
and this finishes the proof statement(ii). \\

(iii) If $M$ is a simple $H$-module with character $\chi$ such that $\chi(\nu_H)\ne 0$, by (ii),
$M$ admits a non-degenerate bilinear
form $\<\cdot, \cdot\>$ with adjoint $S$ such that
$$
\<x,y\> = \chi(\nu_H) \<y, g^{-1}x \>
$$
for all $x, y \in M$. Thus, we have
$$
\begin{aligned}
\<x, y\> & = \chi(\nu_H)^2 \< g^{-1}x, g^{-1}y\>= \chi(\nu_H)^2
\<x, S(g^{-1})g^{-1}y\>\\
& = \chi(\nu_H)^2 \<x, y\> \quad (\mbox{ by  Theorem }\ref{t:7.2} )\,.
\end{aligned}
$$
Therefore,  $\chi(\nu_H)^2=1$ or equivalently $\chi(\nu_H)=\pm 1$.\\

Let $\sum_i x_i \otimes y_i= q_R \Delta(\Lambda)p_R$ where $\Lambda$ is the
normalized two-sided integral of $H$. By Lemma \ref{l:dualbasis},
$
\{S(x_i)g^{-1}, y_i\}
$
is a pair of dual bases of $H$ with respect to the
form $\<\cdot , \cdot \>_{reg}$ on $H$. Therefore, we obtain
\begin{eqnarray*}
\Tr(S) & =& \sum_i \<S(S(x_i)g^{-1}), y_i\>_{reg} \\
&=& \sum_i \chi_{reg}(S(g^{-1})S^2(x_i)y_i)\\
&=& \sum_i\chi_{reg}(S(g^{-1})g^{-1} x_i g y_i)\\
&=& \sum_i \chi_{reg}(S(g^{-1})g^{-1} S(g)  x_i  y_i)
\quad\mbox{by Lemma \ref{pq_1}}\\
&=& \sum_i \chi_{reg}(g^{-1}   x_i  y_i)
\quad\mbox{by Lemma \ref{l:inner}}\\
&=& \chi_{reg}(g^{-1}  \nu_H)\,.
\end{eqnarray*}
Since $\nu_H$ is in the center of $H$, for any
irreducible $H$-module $M$ with character $\chi$, $\nu_H$ acts
on $M$ as a multiplication by the scalar
$$
c_\chi = \chi(\nu_H)/\chi(1_H)\,.
$$
Since $\displaystyle \chi_{reg}=\sum_{\chi \in Irr(H)} \chi(1_H)\chi$, we have
\begin{eqnarray*}
\Tr(S) & =  &\sum_{\chi \in Irr(H)}\chi(1_H)
\chi(g^{-1}   \nu_H) \\
&=& \sum_{\chi \in Irr(H)}\chi(1_H)
c_\chi \,\chi(g^{-1})\\
&=& \sum_{\chi \in Irr(H)}\chi(\nu_H)\chi(g^{-1})\,.
\end{eqnarray*}
\hfill\qed
\begin{remark}
{\rm
In \cite{FGSV99}, Fuchs et al also define a notion of Frobenius-Schur indictor
for
simple objects in a sovereign $C^*$-category $\mathcal{C}$ such that
$$
id : {^*\!M} \map M^*
$$
defines an isomorphism of the tensor functors ${^*?}$ and ${?^*}$.
Let $k_M : M \map (^*\!M)^*={^{**}\!M}$ be the natural isomorphism of the
underlying autonomous structure of $\mathcal{C}$. Then for any simple object
$M$ in $\mathcal{C}$,  the Frobenius-Schur indicator $c_M$ of $M$ is defined to
be 0 if $M \not\cong {^*M}$ and $c$ if there exists a $H$-module isomorphism
$J: M \longrightarrow {^*\!M}$ where $c$ given by the equation
\begin{equation}\label{e:Fuch_FS}
J^* \circ k_M = c \, J \,,
\end{equation}
in which case the values of $c_M$ can only be 0, 1 or -1.\\

The category $H$-$\mod_{fin}$ is not of this kind in general. Nevertheless,
if one replaces $k_M$ in (\ref{e:Fuch_FS}) by $j_M : M \map\, {^{**}\!M}$,
 given by
$$
j_M(f)(x)=f(g^{-1}x) \quad \mbox{for all } x\in M\, \mbox{ and }f \in {^*\!M}\,,
$$
one can still define
{\em Frobenius-Schur
indicator} $c_M$ for any simple $H$-module $M$ to be 0 if
$M \not\cong {^*M}$ and $c$ if there exists a $H$-module isomorphism
$J: M \longrightarrow {^*M}$ where $c$ given by the equation
$$
J^* \circ j_M = c\, J \,.
$$
Theorem \ref{t:main2} (i) and (ii)  implies $c_M = \chi(\nu_H)$. \qed\\
}
\end{remark}

Before closing this section, we will  show that
if $\a$ is a central unit, a bilinear form  on a
$H$-module $M$ is $H$-invariant if, and only if, $S$ is the adjoint
of the form. Both  semi-simple Hopf algebras over $k$ or
 twisted quantum doubles of finite groups are of this type.

\begin{prop} \label{p:adj_inv}
Let $H=(H, \Delta,\e, \Phi, \a, \b, S)$ be a quasi-Hopf algebra
over $k$ and $M$ a $H$-module. Then, the set $Inv(M)$ of $H$-invariant forms
on $M$ and the set $Adj_S(M)$ of  forms $M$ with adjoint $S$ are isomorphic as
$k$-spaces. In addition, if  $\a$ is a central unit of $H$, then
$$
Inv(M) = Adj_S(M)\,.
$$
\end{prop}
\pf Note that both $Inv(M)$ and $Adj_S(M)$ are $k$-subspaces of
$(M \otimes M)^*$. We define $\phi: Adj_S(M) \map (M \otimes M)^*$ and
$\psi:  Inv(M) \map (M \otimes M)^*$ by
\begin{eqnarray}
\phi(\mathbf{b})(x \otimes y) &=& \mathbf{b}(x \otimes \a y)\\
\psi(\mathbf{b'})(x \otimes y) &=& \mathbf{b'}(p_L(x \otimes y))
\end{eqnarray}
for any $x, y \in M$, $\mathbf{b} \in Adj_S(M)$ and
$\mathbf{b'}\in Inv(M)$. Using (\ref{2.8}), one can easily see that
$$
Im (\phi) \C Inv(M)\,.
$$
By (\ref{eq:pL}),
$\psi(\mathbf{b'})$ has adjoint $S$ for any
$H$-invariant form $\mathbf{b'}$ on $M$ and so
$$
Im (\psi) \C Adj_S(M)\,.
$$

It follows easily from (\ref{2.18}) that
for any $\mathbf{b'} \in Inv(M)$ and $x, y \in M$,
\begin{eqnarray*}
\mathbf{b'}(x\otimes y)& =&
\mathbf{b'}(\Delta(q_L^2) p_L (S^{-1}(q_L^1) x \otimes y)) \\
& =&  \mathbf{b'}(p_L (S^{-1}(q_L^1 \e(q_L^2)) x \otimes y) \\
&=& \psi(\mathbf{b'})((S^{-1}(\a) x \otimes y)\,.
\end{eqnarray*}
Since $\psi(\mathbf{b'}) \in Adj_S(M)$,
$$
\phi\circ\psi = id_{Inv(M)}\,.
$$
On the other hand, by (\ref{2.9}), for any
$\mathbf{b} \in Adj_S(M)$ and $x, y \in M$,
$$
\psi\circ \phi(\mathbf{b})(x \otimes y) =
\mathbf{b}(p^1_L x \otimes \a p^2_L  y) =
\mathbf{b}(x \otimes S(p_L^1)\a p^2_L  y) = \mathbf{b}(x \otimes y)\,.
$$
Therefore, $\phi: Adj_S(M) \map Inv(M)$ is a $k$-linear isomorphism. \\

If $\a$ is a central unit, we consider the  quasi-Hopf algebra
$H_{\a^{-1}}$. Then, the corresponding $\phi$ is the identity map and so
$$
Adj_{S_{\a^{-1}}}(M)=Inv(M)\,.
$$
Since $S_{\a^{-1}} =S$, the second statement follows. \qed

\section{Frobenius-Schur Indicators of Twisted Quantum Doubles
of Finite Groups} We showed in section \ref{s:Bantay_Indicator}
that for any simple module $M$ for $D^\w(G)$ with character
$\chi$,  Bantay's formula of the indicator of
$\chi$  is $\chi(\nu_{D^\w(G)})$.   In this section, we will prove that
the trace element of $D^\w(G)$ is $\b$ and
the Frobenius-Schur indictor $\chi(\nu_{D^{\w}(G)})$ of $\chi$ is
non-zero if, and only if, ${^*\!M} \cong M$.
Moreover, the indicator of $\chi$ is
$1$ (respectively $-1$) if and only if $M$ admits a $\b^{-1}$-symmetric
(resp. $\b^{-1}$-skew symmetric) non-degenerate $D^{\w}(G)$-invariant
bilinear form $\<\cdot, \cdot\>$, that is
$$
\<x,y\>= \<y, \b^{-1} x\>\quad (resp. \quad \<x,y\>= -\<y, \b^{-1} x\> )
$$
for all $x, y \in M$.\\

We first need the following formula (cf. \cite{AC92})
to compute the trace element of $D^\w(G)$.
\begin{lem}\label{l:8.1}
Let $\w:G \times G \times G \map k^\times $ be a normalized
3-cocycle of a finite group $G$ and let $S$ be the antipode of
the quasi-Hopf algebra $D^\w(G)$ defined in Section \ref{s:Bantay_Indicator}.
Then for any $g, x \in G$,
\begin{eqnarray*}
S^2(e(g)\otimes x) & = & \frac{\w(g^x, (g^{-1})^x, g^x)}{\w(g, g^{-1},g)}
\, e(g)\otimes x\,,\\
&=& \b^{-1} (e(g) \otimes x) \b\,.
\end{eqnarray*}
\end{lem}
\pf It follows from (\ref{eq:antipode}) that
\begin{equation}\label{eq:S2}
S^2(e(g)\otimes x) = \left(
\theta_{g^{-1}}(x,x^{-1})\g_x(g,g^{-1})\theta_{g^x}(x^{-1},
x)\g_{x^{-1}}((g^{-1})^x,g^x) \right)^{-1} e(g)\otimes x\,,
\end{equation}
where $g^x$ denotes the product $x^{-1}gx$.  By the
normality of $\w$ and (\ref{eq:0.1}),
$$
\theta_g(x, x^{-1}) = \theta_{g^x}(x^{-1}, x)\,.
$$
 Thus, we
have
\begin{eqnarray*}
&&\theta_{g^{-1}}(x,x^{-1})\g_x(g,g^{-1})\theta_{g^x}(x^{-1}, x)
  \g_{x^{-1}}((g^{-1})^x,g^x) \\ \\
& = &\theta_{g^{-1}}(x,x^{-1})\theta_{g}(x, x^{-1})
   \g_x(g,g^{-1})\g_{x^{-1}}((g^{-1})^x,g^x)
\end{eqnarray*}
By the normality of $\w$ and equation
(\ref{eq:0.2}), we have
\begin{eqnarray*}
&&\theta_{g^{-1}}(x,x^{-1})\g_x(g,g^{-1})\theta_{g^x}(x^{-1}, x)
  \g_{x^{-1}}((g^{-1})^x,g^x)\\ \\
& = & \frac{\g_x(g,g^{-1})\g_{x^{-1}}((g^{-1})^x,g^x)}{\g_x(g, g^{-1})\g_{x^{-1}}(g^x, (g^{-1})^x)}\\ \\
&=& \frac{\g_{x^{-1}}((g^{-1})^x,g^x)}{\g_{x^{-1}}(g^x, (g^{-1})^x)}\,.
\end{eqnarray*}
By equation (\ref{eq:0.3}), for any $z, a \in G$ we have
$$
\g_z(a,a^{-1}) \w(a^z, (a^{-1})^z, a^z)=\g_z(a^{-1}, a)\w(a,a^{-1},a)\,.
$$
Hence we have
\begin{eqnarray}\nonumber
\frac{\g_{x^{-1}}((g^{-1})^x,g^x)}{\g_{x^{-1}}(g^x, (g^{-1})^x)} &=&
\frac{\w((g^x)^{x^{-1}}, ((g^{-1})^x)^{x^{-1}}, (g^x)^{x^{-1}})}{\w(g^x, (g^{-1})^x, g^x)}\\
\label{eq:0.4}  \\\nonumber
&=&
\frac{\w(g, g^{-1}, g)}{\w(g^x, (g^{-1})^x, g^x)}\,.\\ \nonumber
\end{eqnarray}
The second equation in the statement of the Lemma follows
immediately from
(\ref{multiplication}). \qed\\

\begin{prop}
Let $\w:G \times G \times G \map k^\times $ be a normalized
3-cocycle of a finite group $G$. Then the trace element of the
quasi-Hopf algebra $D^\w(G)$ is $\b$.
\end{prop}
\pf Using (\ref{eq:antipode}), $S(\b)=\b^{-1}$. Suppose that
$\{f_{g,x}\}_{g,x \in G}$ is the dual basis of $\{e(g) \otimes
x\}_{g,x \in G}$. Then, by (\ref{eq:cointegral}), the normalized
left cointegral of $D^\w(G)$ is  given by
$$
 \lambda(e(g) \otimes x)  =  \sum_{h,y \in G} f_{h,y}((e(g) \otimes x)
S^2(e(h) \otimes y)\b^{-1} )\,.
$$
Using Lemma \ref{l:8.1}, we have
\begin{eqnarray*}
 \lambda(e(g) \otimes x) & = & \sum_{h,y \in G} f_{h,y}((e(g) \otimes x)
 \b^{-1} (e(h) \otimes y ))\\
 & =& \chi_{reg}((e(g) \otimes x)\b^{-1})\\
 & = & \lambda((e(g) \otimes x)\b^{-1} g) \,.
\end{eqnarray*}
By the non-degeneracy of $\lambda$, $\b^{-1} g =1$ and so $g= \b$. \qed\\

\begin{cor}
Let $\w:G \times G \times G \map k^\times $ be a normalized
3-cocycle of a finite group $G$. Suppose that $M$ is
a simple $D^\w(G)$-module with character $\chi$. Then the
Frobenius-Schur indicator $\chi(\nu_{D^\w(G)})$ of $\chi$ satisfies
the following properties:
\begin{enumerate}
\item[\rm (i)] $\chi(\nu_{D^\w(G)}) = 0, 1$, or $-1$.
\item[\rm (ii)] $ \chi(\nu_{D^\w(G)}) \neq 0$ if, and only if,
${^*\!M} \cong M$\,.
\item[\rm (iii)] $\chi(\nu_{D^\w(G)})=1$ (respectively $-1$) if
and only if $M$ admits a $\b^{-1}$-symmetric (resp. $\b^{-1}$-skew symmetric)
non-degenerate $D^{\omega}(G)$-invariant bilinear form.
\end{enumerate}
Moreover,
$$
\Tr(S) = \sum_{\chi \in Irr(D^\w(G))} \chi(\nu_{D^\w(G)}) \chi(\b^{-1})\,.
$$
\end{cor}
\pf Statement (i), (ii) and the last
statement are immediate consequences of Theorem
\ref{t:main2}. Since $\a=1$, by Proposition \ref{p:adj_inv},
a bilinear form $\<\cdot, \cdot\>$ on $M$ is $D^\w(G)$-invariant if,
and only if, $S$ is the adjoint of $\<\cdot, \cdot\>$. Thus, by
Theorem \ref{t:main2} (ii), the result in statement (iii) follows. \qed

%GATHER{MYBIBL.BIB}
\bibliographystyle{amsalpha}
%\bibliography{mybibl}

\begin{thebibliography}{FGSV99}

\bibitem[AC92]{AC92}
Daniel Altsch{\"u}ler and Antoine Coste, \emph{Quasi-quantum groups, knots,
  three-manifolds, and topological field theory}, Comm. Math. Phys.
  \textbf{150} (1992), no.~1, 83--107. \MR{94b:57006}

\bibitem[Ban97]{Bantay97}
Peter Bantay, \emph{The {F}robenius-{S}chur indicator in conformal field
  theory}, Phys. Lett. B \textbf{394} (1997), no.~1-2, 87--88. \MR{98c:81195}

\bibitem[Ban00]{Bantay00}
\bysame, \emph{Frobenius-{S}chur indicators, the {K}lein-bottle amplitude, and
  the principle of orbifold covariance}, Phys. Lett. B \textbf{488} (2000),
  no.~2, 207--210. \MR{2001e:81094}

\bibitem[Ban02]{BantayConv}
\bysame, \emph{Private communication, 2002}.

\bibitem[BK01]{BaKi}
Bojko Bakalov and Alexander Kirillov, Jr., \emph{Lectures on tensor categories
  and modular functors}, University Lecture Series, vol.~21, American
  Mathematical Society, Providence, RI, 2001. \MR{2002d:18003}

\bibitem[CR88]{CR88}
Charles~W. Curtis and Irving Reiner, \emph{Representation theory of finite
  groups and associative algebras}, Wiley Classics Library, John Wiley \& Sons
  Inc., New York, 1988, Reprint of the 1962 original, A Wiley-Interscience
  Publication. \MR{90g:16001}

\bibitem[DPR92]{DPR90}
R.~Dijkgraaf, V.~Pasquier, and P.~Roche, \emph{Quasi-{H}opf algebras, group
  cohomology and orbifold models [{M}{R} 92m:81238]}, Integrable systems and
  quantum groups (Pavia, 1990), World Sci. Publishing, River Edge, NJ, 1992,
  pp.~75--98.

\bibitem[Dri90]{Drin90}
V.~G. Drinfel'd, \emph{Quasi-{H}opf algebras}, Leningrad Math. J. \textbf{1}
  (1990), 1419--1457.

\bibitem[EG02]{EG02}
Pavel Etingof and Shlomo Gelaki, \emph{On families of triangular {H}opf
  algebras}, Int. Math. Res. Not. (2002), no.~14, 757--768. \MR{2002m:16036}

\bibitem[ENO]{ENO}
Pavel Etingof, Dmitri Nikshych, and Viktor Ostrik, \emph{On fusion categories},
  preprint \textbf{arXiv:math.QA/0203060}.

\bibitem[FGSV99]{FGSV99}
J.~Fuchs, A.~Ch. Ganchev, K.~Szlach{\'a}nyi, and P.~Vecserny{\'e}s,
  \emph{{$S\sb 4$} symmetry of {$6j$} symbols and {F}robenius-{S}chur
  indicators in rigid monoidal {$C\sp *$} categories}, J. Math. Phys.
  \textbf{40} (1999), no.~1, 408--426. \MR{99k:81111}

\bibitem[FY92]{FY}
Peter Freyd and David~N. Yetter, \emph{Coherence theorems via knot theory}, J.
  Pure Appl. Algebra \textbf{78} (1992), no.~1, 49--76. \MR{93d:18013}

\bibitem[HN99a]{HN992}
Frank Hausser and Florian Nill, \emph{Diagonal crossed products by duals of
  quasi-quantum groups}, Rev. Math. Phys. \textbf{11} (1999), no.~5, 553--629.
  \MR{2000d:81069}

\bibitem[HN99b]{HN991}
\bysame, \emph{Doubles of quasi-quantum groups}, Comm. Math. Phys. \textbf{199}
  (1999), no.~3, 547--589. \MR{2000a:16075}

\bibitem[HN]{HNQA}
\bysame, \emph{Integral theory for quasi-hopf algebras}, preprint
  \textbf{arXiv.math.QA/9904164}.

\bibitem[Kas95]{Kassel}
Christian Kassel, \emph{Quantum groups}, Springer-Verlag, New York, 1995.

\bibitem[KMM02]{KMM02}
Y.~Kashina, G.~Mason, and S.~Montgomery, \emph{Computing the
  {F}robenius-{S}chur indicator for abelian extensions of {H}opf algebras}, J.
  Algebra \textbf{251} (2002), no.~2, 888--913. \MR{1 919 158}

\bibitem[LM00]{LM00}
V.~Linchenko and S.~Montgomery, \emph{A {F}robenius-{S}chur theorem for {H}opf
  algebras}, Algebr. Represent. Theory \textbf{3} (2000), no.~4, 347--355,
  Special issue dedicated to Klaus Roggenkamp on the occasion of his 60th
  birthday. \MR{2001k:16073}

\bibitem[LR87]{LaRa87}
Richard~G. Larson and David~E. Radford, \emph{Semisimple cosemisimple {H}opf
  algebras}, Amer. J. Math. \textbf{109} (1987), no.~1, 187--195.
  \MR{89a:16011}

\bibitem[LR88]{LaRa88}
\bysame, \emph{Finite-dimensional cosemisimple {H}opf algebras in
  characteristic $0$ are semisimple}, J. Algebra \textbf{117} (1988), no.~2,
  267--289. \MR{89k:16016}

\bibitem[Mas95]{Mas}
Geoffrey Mason, \emph{The quantum double of a finite group and its role in
  conformal field theory}, Groups '93 Galway/St.\ Andrews, Vol.\ 2, Cambridge
  Univ. Press, Cambridge, 1995, pp.~405--417. \MR{97a:11067}

\bibitem[MN01]{MN01}
Geoffrey Mason and Siu-Hung Ng, \emph{Group cohomology and gauge equivalence of
  some twisted quantum doubles}, Trans. Amer. Math. Soc. \textbf{353} (2001),
  no.~9, 3465--3509 (electronic). \MR{1 837 244}

\bibitem[MN]{MN03}
\bysame, \emph{Central invariants and frobenius-schur indicators for semisimple
  quasi-hopf algebras}, preprint \textbf{arXiv.math.QA/0304156}.

\bibitem[Mue]{Mu}
Michael Mueger, \emph{From subfactors to categories and topology i. frobenius
  algebras in and morita equivalence of tensor categories}, preprint
  \textbf{arXiv:math.CT/0111204}.

\bibitem[Pan98]{Pana98}
Florin Panaite, \emph{A {M}aschke-type theorem for quasi-{H}opf algebras},
  Rings, Hopf algebras, and Brauer groups (Antwerp/Brussels, 1996), Lecture
  Notes in Pure and Appl. Math., vol. 197, Dekker, New York, 1998,
  pp.~201--207. \MR{99k:16085}

\bibitem[PVO00]{PV00}
Florin Panaite and Freddy Van~Oystaeyen, \emph{Existence of integrals for
  finite dimensional quasi-{H}opf algebras}, Bull. Belg. Math. Soc. Simon
  Stevin \textbf{7} (2000), no.~2, 261--264. \MR{2001f:16079}

\bibitem[Ser77]{serre77}
Jean-Pierre Serre, \emph{Linear representations of finite groups},
  Springer-Verlag, New York, 1977, Translated from the second French edition by
  Leonard L. Scott, Graduate Texts in Mathematics, Vol. 42. \MR{56 \#8675}

\bibitem[TY98]{TaYa98}
Daisuke Tambara and Shigeru Yamagami, \emph{Tensor categories with fusion rules
  of self-duality for finite abelian groups}, J. Algebra \textbf{209} (1998),
  no.~2, 692--707.

\end{thebibliography}
\providecommand{\bysame}{\leavevmode\hbox to3em{\hrulefill}\thinspace}
\providecommand{\MR}{\relax\ifhmode\unskip\space\fi MR }
\providecommand{\MRhref}[2]{%
  \href{http://www.ams.org/mathscinet-getitem?mr=#1}{#2}
}
\providecommand{\href}[2]{#2}

\end{document}